\documentclass[aos, preprint]{imsart}
\usepackage[utf8]{inputenc}
\usepackage{amsmath, amssymb, amsthm, amsfonts, comment, mathrsfs, graphicx, hyperref, dsfont, tikz}
\usepackage[nothm]{yk}
\usepackage[authoryear]{natbib}
\RequirePackage[colorlinks,citecolor=blue,urlcolor=blue]{hyperref}
\usepackage{neuralnetwork}
\newcommand{\fn}{\mathrm{final}}
\usepackage{bbm, framed}
\usepackage[inline]{enumitem}   
\usepackage{xcolor}
\usepackage[linesnumbered,ruled,vlined]{algorithm2e}
\SetKwInput{KwInput}{Input}                
\SetKwInput{KwOutput}{Output}              

\startlocaldefs
\theoremstyle{plain}

\newtheorem{theorem}{Theorem}[section]
\newtheorem{lemma}[theorem]{Lemma}
\newtheorem{assumption}[theorem]{Assumption}
\newtheorem{remark}[theorem]{Remark}
\newtheorem{corollary}[theorem]{Corollary}
\newtheorem{proposition}[theorem]{Proposition}
\theoremstyle{remark}

\endlocaldefs

\begin{document}
\begin{frontmatter}

\title{Deep Neural Networks for Nonparametric Interaction Models with Diverging Dimension}

\runtitle{DNN for Nonparametric Interaction Models}
\begin{aug}
\author[A]{\fnms{Sohom}~\snm{Bhattacharya}\ead[label=e1]{}},
\author[A]{\fnms{Jianqing}~\snm{Fan}\ead[label=e1]{}}
\and
\author[A]{\fnms{Debarghya}~\snm{Mukherjee}\ead[label=e1]{}}
\address[A]{Department of Operations Research and Financial Engineering,
Princeton University\printead[presep={,\ }]{e1}}
\end{aug}
\date{}

\begin{abstract}
\noindent Deep neural networks have achieved tremendous success due to their representation power and adaptation to low-dimensional structures.  Their potential for estimating structured regression functions has been recently established in the literature. However, most of the studies require the input dimension to be fixed and consequently ignore the effect of dimension on the rate of convergence and hamper their applications to modern big data with high dimensionality. In this paper, we bridge this gap by analyzing a $k^{th}$ order nonparametric interaction model in both growing dimension scenarios ($d$ grows with $n$ but at a slower rate) and in high dimension ($d \gtrsim n$).  In the latter case, sparsity assumptions and associated regularization are required in order to obtain optimal rates of convergence.
A new challenge in diverging dimension setting is in calculation mean-square error, the covariance terms among estimated additive components are an order of magnitude larger than those of the variances and they can deteriorate statistical properties without proper care.  We introduce a critical debiasing technique to amend the problem. We show that under certain standard assumptions, debiased deep neural networks achieve a minimax optimal rate both in terms of $(n, d)$. Our proof techniques rely crucially on a novel debiasing technique that makes the covariances of additive components negligible in the mean-square error calculation.  In addition, we establish the matching lower bounds. 
\end{abstract}



\end{frontmatter}

\section{Introduction}

Recent advances in technology have allowed
statisticians to collect data on a large number of explanatory variables to
predict outcomes of interest~\citep{goodfellow2016deep,fan2020statistical}. Often, the relationship between these outcomes and predictors are highly nonlinear (e.g., image data like MNIST \citep{lecun1998mnist}, CIFAR \citep{krizhevsky2009learning} etc.), yielding a practical need to investigate multivariate nonparametric regression model
\begin{equation}
    Y_i= f(X_i)+ \varepsilon_i, \quad i=1,\ldots,n,
\end{equation}
where $X_i$'s are explanatory variables, $Y_i$'s are response variables and $\varepsilon_i$'s are unobserved errors. The statistical problem here is to recover $f$ under some minimal smoothness assumptions. A classical result of Charles Stone~\citep{stone1982optimal} shows if $f$ is a $d$-variate $(\beta,C)$ smooth function (precise definition given later), the minimax optimal rate of estimation is of the order $n^{-\frac{2\beta}{2\beta+d}}$, which is referred to as ``curse of dimensionality", i.e., large sample size is necessary to estimate the regression function well. In particular, when $\beta$ is fixed, it is easily to see that when $d \asymp \log n$, minimax consistent estimators can not be obtained.

To tackle aforementioned issue, one needs to impose a low-dimensional structure on the function $f$~\citep{kpotufe2012tree,yang2016bayesian, yang2015minimax}, such as additive regression model \citep{stone1985additive}, projection pursuit regression model \citep{friedman1981projection}, higher order interaction model \citep{stone1994use}, generalized higher-order interaction model \citep{horowitz2007rate, kohler2016nonparametric}, single and multi-index models \citep{hardle1993optimal}, etc.; see \cite{fan1996local} for an overview. The key idea behind all these simplifications is to reduce the complexity of the underlying function class and mollify the effect of dimension in the rate of convergence. For example, \cite{stone1985additive} proved that for additive regression model (when $f(x) = \sum_{j=1}^d f_j(x_j)$) with $\beta$-smooth univariate components, the minimax optimal rate of estimation of $f$ in terms of squared $L_2(P_X)$ loss is $C_dn^{-\frac{2\beta}{2\beta + 1}}$. \cite{fan1998direct} demonstrated further the adaptivity of such an additive structure and \cite{horowitz2007rate} studied further a general class of nonparametric regression models with unknown link functions. As each component function $f_j$ is univariate and can be estimated at a rate $n^{-\frac{2\beta}{2\beta + 1}}$, the effect of dimension appears only through the multiplicative constant $C_d$, not in the power of $n$. Later, \cite{stone1994use} extended this result for higher-order interaction model (i.e. when $f(x) = \sum_{J \in S} f_J(x_J)$, $S$ is a collection of subsets of $\{1, \dots, d\}$ and $|J| \le d^\star$ for all $J \in S$) and showed that  the minimax optimal rate of estimation is $C_dn^{-\frac{2\beta}{2\beta + d^\star}}$, where, as before, the effect of $d$ appears through a multiplicative constant. Again, all of the above results assume that $d$ is finite, which is inappropriate in many modern data science applications.

In practice, it is not enough just to obtain minimax optimal rates, efficient and easily computable estimators are warranted. Several methods like kernel, spline, and wavelet based estimators (see Section 2 of \cite{fan1996local} for an overview), backfitting algorithm \citep{buja1989linear},  off-the-shelf nonparametric methods like boosting \citep{freund1996experiments}, random forest \citep{breiman2001random} have been popularly used. Recently, neural networks have emerged as an imperative tool for analyzing nonlinear relations between covariates and response variables. Deep (or multilayer) neural networks have been the backbone of the incredible advances in machine learning that have resulted in massive success in reinforcement
learning, robotics, computer
vision, natural language processing, and other statistical prediction tasks~\citep {goodfellow2016deep}. 
Leveraging the availability of large amounts of digitized data, deep learning has enjoyed a plethora of empirical successes. However, most of the successes are purely human engineered (i.e. achieved after tuning lots of hyperparameters), lacking theoretical guarantees. As a consequence, researchers are naturally interested in understanding the statistical properties of such useful and practical estimators. In \cite{kohler2005adaptive}, the authors established that for the standard nonparametric models as well as for the higher-order interaction models, deep neural network (henceforth DNN) based estimators are minimax rate optimal up to logarithmic factor. Later, the results are extended for a more general class of functions (namely generalized hierarchical interaction model, see \cite{bauer2019deep}) in a series of work by\cite{kohler2016nonparametric, bauer2019deep, schmidt2020nonparametric, kohler2021rate, fan2022noise}.

The previous research in understanding the rate of convergence of DNN-based estimators assumes the underlying dimension of covariate $d$ is fixed and consequently ignores the effect of $d$ in the rate of the estimator. 
However, in the era of big data, the dimension of underlying covariates is quite large in many practical problems, often larger than the underlying sample size $n$. For example, in a simple image classification problem, a $28 \times 28$ image lies in $\reals^{784}$. Similarly, in genome-wide association studies, the number of single nucleotide polymorphism (SNP) can be significantly higher than the number of individuals\citep{novembre2008genes}.  In that paper, the authors have analyzed a dataset consisting of $1387$ individuals and $\sim 2\times 10^5$ SNPs. In such examples, 
a sharper analysis quantifying the optimal dependence on $d$ is often necessary. 
 
In this paper, we aim to bridge this gap by analyzing the nonparametric regression model both in growing dimension (when $d$ grows with $n$ but $d = o(n)$) and high dimension (when $d \gtrsim n$) setup. The distinction between these two cases lies in whether or not the sparsity and regularization are required for consistent estimation, in an analogue way to the linear model with growing dimension and sparse linear model in high dimension.  We analyze the following $k$-way interaction model \citep{stone1994use}: 
\begin{equation}\label{eq:define_f}
    Y_i = \sum_{l=1}^k \sum_{\substack{(j_1 < \dots <j_l) \in [d]^l}}f_{j_1,j_2,\ldots,j_l}(X_{i,\{j_1,j_2,\ldots,j_l\}}) + \eps_i \,,
\end{equation}
where $[d]^l$ is the collection of all ordered subsequence of length $l$ of $\{1, \dots, d\}$.  As mentioned before, $k = o(\log n)$ is necessary in order to have a consistent estimator.  Hence, we will take a finite $k$.
Such models include the well-studied additive model~\citep{stone1985additive}, namely, $f(x)=\sum_{j=1}^d f_j(x_j)$ and interaction models~\citep{stone1994use}. When $d \gtrsim n$, it is necessary to impose sparsity assumption even in the high-dimensional linear model. We, therefore, consider a \textit{sparse $k$-way interaction model} as follows: 
\begin{equation}\label{eq:define_spam}
    Y_i = \sum_{l=1}^k \sum_{\substack{(j_1,j_2,\ldots,j_l) \in S_l}}f_{j_1,j_2,\ldots,j_l}(X_{i,\{j_1,j_2,\ldots,j_l\}}) + \eps_i \,,
\end{equation}
where $S_l \subset [d]^l$ with $|S_l| \ll d^l$.
When each univariate function is constrained to be linear, the model \eqref{eq:define_spam} reduces to a sparse parametric interaction model. When $k=1$, the model \eqref{eq:define_spam} becomes a \textit{sparse additive model} which has been well-studied in literature~\citep{lin2006component,koltchinskii2010sparsity,ravikumar2009sparse,tan2019doubly,yuan2016minimax}, building upon the recent developments on penalized linear regression. For example, \cite{tan2019doubly} proposes to doubly penalize each component $f_j$ by its empirical norm and functional semi-norm to estimate the regression function.

When $d \ll n$, we estimate the mean function of \eqref{eq:define_f} via minimizing squared error loss (details can be found in Section \ref{sec:analysis_growing_dim}). However, when $d \gtrsim n$, we need to penalize explicitly to enforce sparsity. In previous works on the sparse additive models, researchers typically impose two penalties: i) one to enforce sparsity (e.g. via the summation of $L_2(\bbP_n)$ norm on the component functions) and ii) another to control the complexity of the component functions (e.g., penalize with respect to RKHS norm if the component functions lie in a RKHS). In reality, implementing such a doubly penalized estimator is often computationally challenging. To overcome such difficulty, we introduce a two-step hard thresholding-based estimator, which is motivated by the seminal work of SURE independent screening \citep{fan2008sure, fan2010sure} and least square estimation after model selection in high dimensional linear regression model \citep{belloni2013least}. The idea can be briefly described as follows: first, we obtain an initial estimator (say $\hat f^{\text{init}}$) of the mean function by only penalizing the summation of empirical $L_2(\bbP_n)$ norm of the component functions (which can be implemented by penalizing the last layer of neural network). Analogously, this can also be thought of as a version of a group lasso penalty, where the neural network corresponding to each component is a group of variables and we penalize the sum of the norm of each group (here we take the $L_2(\bbP_n))$ norm) to enforce a component-wise sparsity. 
Next, we further perform a hard thresholding on the non-zero component of $\hat f^{\text{init}}$ to weed out the \emph{small} non-zero components and estimate the active sets. Finally, we perform an empirical risk minimization by minimizing the squared error loss over only the components selected in the previous step to obtain the final estimator. We prove that under some fairly standard assumptions, the estimator is minimax rate optimal. We now summarize our contribution as follows.

{\bf Contribution: } Our main contribution is a rigorous theoretical analysis of $k$-way interaction model (also known as nonparametric ANOVA model) in both growing and high dimensional setup. To achieve this goal, we also prove several results that may be of independent interest. We summarise the key contributions below:  
\begin{itemize}
\item We analyze the nonparametric $k$-way interaction model using neural network based estimator both when $d = \dim(X)$ increases with $n$ and $d = o(n)$, and when $d \gtrsim n$ with aforementioned regularization. We show that the neural network-based estimator is minimax rate optimal up to a poly-log factor. 

\item We introduce a novel debiasing technique that makes the covariances among additive components negligible, which reduces statistical errors.

\item We prove an approximation result of smooth function using deep neural network (Theorem \ref{thm:approx_growing_d}) via a \emph{novel debiasing technique}, which implies that one can approximate a smooth function using neural networks at the same rate even under constraints (i.e. marginals of a multivariate function are $0$). 

\item We establish the minimax lower bound for estimating the high dimensional $k$-way interaction model, which, to the best of our knowledge, is not present in the literature. 
\end{itemize}

The rest of the paper is organized as follows: In Section \ref{sec:analysis_growing_dim}, we analyze the $k$ way interaction model when $d \ll n $. We divide the entire analysis into three subsections: Subsection \ref{sec:approx_lowdim} bounds approximation error, Subsection \ref{sec:stat_lowdim} analyzes the statistical error, and Subsection \ref{sec:minimax_lowdim} develops the minimax lower bound. Section \ref{sec:high_dim_results} deals with the analysis when $d \gtrsim n$. We broadly divide our analysis into two parts; Subsection \ref{sec:fixed_design_highdim} contains the analysis for fixed designs and Subsection \ref{sec:random_design_highdim} deals with the case of random designs. Furthermore, in Subsection \ref{sec:minimax_high_dim}, we establish the minimax lower bound for the sparse nonparametric interaction models.Section \ref{sec:conclusion} provides conclusional remarks. Finally, Section \ref{sec:proof_approx_growing_d} contains the proof of Theorem \ref{thm:approx_growing_d}. The remaining proofs can be found in the Appendix. 

\section{Analysis of DNN in low diverging dimension}
\label{sec:analysis_growing_dim}
\begin{figure}
    \centering
    \def\layersep{2.5cm}
\begin{tikzpicture}[
   shorten >=1pt,->,
   draw=black!50,
    node distance=\layersep,
    every pin edge/.style={<-,shorten <=1pt},
    neuron/.style={circle,fill=black!25,minimum size=17pt,inner sep=0pt},
    input neuron/.style={neuron, fill=blue!50},
    output neuron/.style={neuron, fill=blue!50},
    hidden neuron/.style={neuron, fill=green!50},
    annot/.style={text width=4em, text centered}
]

\foreach \name / \y in {1,...,1}
    \node[input neuron, pin=left:Input] (I-\name) at (0,-1.5) {};

\newcommand\Nhidden{3}

\foreach \N in {1,...,\Nhidden} {
   \foreach \y in {1,...,3} {
      \path[yshift=0.5cm]
          node[hidden neuron] (H\N-\y) at (\N*\layersep,-\y cm) {};
       }
\node[annot,above of=H\N-1, node distance=1cm] (hl\N) {Hidden layer \N};
}

\foreach \name / \y in {1,...,1}
    \node[output neuron, pin=right:Output] (O-\name) at (4*\layersep,-1.5) {};

\foreach \source in {1,...,1}
    \foreach \dest in {1,...,3}
        \path (I-\source) edge (H1-\dest);

\foreach [remember=\N as \lastN (initially 1)] \N in {2,...,\Nhidden}
   \foreach \source in {1,...,3}
       \foreach \dest in {1,...,3}
           \path (H\lastN-\source) edge (H\N-\dest);

\foreach [remember=\N as \lastN (initially 3)] \N in {\Nhidden,...,\Nhidden}
\foreach \source in {1,...,3}
    \foreach \dest in {1,...,1}
        \path (H\lastN-\source) edge (O-\dest);


\node[annot,left of=hl1] {Input layer};
\node[annot,right of=hl\Nhidden] {Output layer};
\end{tikzpicture}
    \caption{A graphic representation of a deep neural network with 3 hidden layers, one input layer, one output layer.}
    \label{fig:deepnet_generic}
\end{figure}
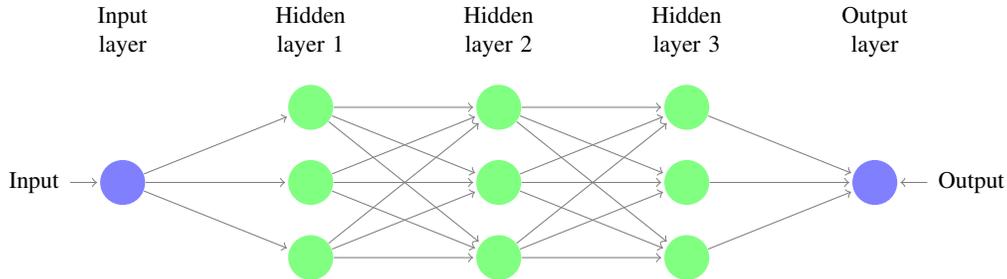

In this section, we present our analysis of the DNN-based estimator of the mean function when the dimension $d$ of $X$ grows with $n$ but $d \ll n$. We consider a $k$-order interaction nonparametric regression model. For an input-output pair $(X, Y)$, where $X \in \reals^d$ and $Y \in \reals$, a $k$-order interaction model is defined as:   
\begin{equation}
\label{eq:model_k}
Y_i = \sum_{\substack{J \subset [d] \\ |J| \le k}}f_{0, J}(X_{i,J}) + \eps_i \,.
\end{equation}
 Here, $\eps_i$'s are assumed to be centered error independent of $X$ and the functions $f_{0, J}: \reals^{|J|} \mapsto \reals$. For simplicity of exposition, we henceforth confine ourselves to a $2$-order interaction model. The extension of our analysis from a $2$-order interaction model to a general $k$-order interaction model is purely technical and the proof ideas will be outlined at the end of the relevant sections. For $k = 2$, the the model presented in equation \eqref{eq:model_k} simplifies to: 
\begin{equation}
\label{eq:two_way_ss_anova}
Y_i = \sum_{j=1}^d f_{0, j}(X_{i,j}) + \sum_{j < k}f_{0, jk}(X_{i, j}, X_{i, k}) + \eps_i \equiv f_0(X_i) + \varepsilon_i \,.
\end{equation}  
where $f_{0, j}'s$ are univariate functions and $f_{0, jk}$'s are bivariate functions. We will make the following assumption on $X$ and $\varepsilon$:

\begin{assumption}
\label{assn:density_bound}
$X$ is supported on $[0, 1]^d$ and admits a density function $p$ such that $\sup_{x \in [0,1]^d}p(x)=: p_{\max} < \infty$, where $p_{\max}$ is free of $d$. We assume $\eps$ is sub-gaussian with sub-gaussian constant $\sigma^2_\eps$. 
\end{assumption}
\noindent Assumption \ref	{assn:density_bound} is a standard assumption in the literature of nonparametric regression. Note that assuming $X$ has compact support is equivalent to assuming that the support is $[0, 1]^d$ via centering and scaling. Many of our results can be extended to the case of unbounded $X$ via truncation arguments -- we omit such arguments for the sake of the simplicity of our article. The upper bound assumption on the density is quite natural, as this merely rules out certain degenerate distributions (i.e. when the density of some region diverges to infinity).  The assumption of sub-gaussianity on the error terms is natural as we aim to analyze the least square estimate of $f$ based on neural networks. It is well-established in the literature that for heavy-tailed errors, a least square estimator of a nonparametric function is not minimax optimal (e.g., see \cite{han2019convergence}) and one should use a variant of Huber loss function. We believe most of our analysis can be extended to this heavy-tailed error, but this is out of the scope of the current paper. The recent article \citep{fan2022noise} sheds some light on the estimation error of deep neural networks with heavy-tailed error in fixed dimension setup. 
\\\\
In our model \eqref{eq:two_way_ss_anova}, the component functions of $f_0$ are not identifiable without further assumptions for two reasons: 
\begin{enumerate}
\item All the functions are identifiable up to shift, i.e. we cannot identify the difference between $(f_{0, i_1 i_2}, f_{0, j_1 j_2})$ and $(f_{0, i_1 i_2} +c, f_{0, j_1 j_2} - c)$ for a constant $c$. 
\item The univariate marginals of the bivariate functions are not identifiable. For example, consider the subcollection of functions $\cC_1 = \{f_{0, 1}, f_{0, 12}, f_{0, 13}, \dots, f_{0, 1d}\}$. Let $g_k(x) = \int f_{0, 1k}(x, y) \ dy$ be the marginal of $f_{0, 1k}$. Then we cannot differentiate between $\cC_1$ and a modified collection $\cC_2 = \{f_{0, 1} + \sum_{k \geq 2}g_k, f_{0, 12} - g_2, \dots, f_{0, 1d} - g_d\}$. 
\end{enumerate}
Therefore, to identify and estimate the non-parametric functions, we need to impose certain structural conditions which will prevent us from shifting. Toward that end, we assume the following: 
\begin{assumption}[Identifiability and boundedness]
\label{assn:identifiability}
We assume the following conditions for identifiability purpose and boundedness: 
\begin{enumerate}
\item All the univariate functions in \eqref{eq:two_way_ss_anova} integrates to $0$, i.e. $\int_0^1 f_{0, j}(x) \ dx = 0$ for $1 \le j \le d$. 
\item All the bivariate functions have $0$ marginals, i.e. 
$$
\int_0^1 f_{0, ij}(x, y) \ dx = \int^1_0 f_{0, ij}(x, y) \ dy = 0 \,.
$$
\item $\|f_0\|_\infty \le B$ for some $B > 0$, where $f_0$ is defined in \eqref{eq:two_way_ss_anova}.
\end{enumerate}
\end{assumption}
\noindent Note that, the above assumption implies that the total integral of $f_{0, ij}$ is also $0$. Assumption \ref{assn:identifiability} ensures identifiability of all the functions involved in \eqref{eq:two_way_ss_anova}. Our next assumption is a smoothness assumption on the underlying component function: 
\begin{assumption}
\label{assn:smoothness}
The functions $\{f_{0, j}\}$ and $\{f_{0, ij}\}$ are assumed to be $(\beta, L)$-Holder smooth, i.e. the functions are $\lfloor \beta 
\rfloor$ times differentiable and $\lfloor \beta \rfloor^{th}$ derivative are Holder with exponent $\beta - \lfloor \beta \rfloor$ and constant $L$. 
\end{assumption}
\noindent The smoothness assumption is standard in the nonparametric regression literature, cf.~\cite[Chapter 1]{tsybakov2004introduction}, as this smoothness assumption controls the complexity of the underlying function class. Another common assumption is that all $f_{0,i}$ and $f_{0,ij}$ belong to a Reproducing Kernel Hilbert Space(RKHS), e.g, see \cite{raskutti2012minimax}, \cite{koltchinskii2010sparsity} and references therein. It would be interesting to study the analog of our results under such assumptions.
\begin{remark}
    Assuming that $\{f_{0, ij}\}$ have a different level of smoothness $\beta_{ij}$ is equivalent to assuming they all have $\min_{i,j} \beta_{i,j}$ smoothness as long as one is concerned about estimation error bounds. 
\end{remark}
\noindent To estimate the mean function of \eqref{eq:two_way_ss_anova} using deep neural networks, we need some properties of DNNs. A neural network is a parametric function $f_\theta$ which maps the input space to output space, i.e. in our context $f_\theta$ maps $\reals^d$ to $\reals$. A two-layer (i.e. one hidden layer) neural network (often termed as \emph{shallow neural network}) with $N$ neurons and activation function $\sigma$ is defined as: 
$$
f_\theta(x) = a + \sum_{i=1}^N b_i \sigma(c_i^\top x + d_i) \,. 
$$
where $\theta = \{a, b_i, c_i, d_i\}$. Therefore, $f_\theta$ first projects $x$ into $\reals^N$ via a linear transformation $x \mapsto Cx + d$ where $c_i^\top$'s are the rows of $C$ and $d_i$'s are elements of $d$. Then it applies a non-linear activation $\sigma$ to all the coordinates of $Cx + d$. Finally, it takes a linear combination of the coordinates using the map $x \mapsto a + \langle b, x 
\rangle$. This shallow neural network can be extended to the deep neural network by adding more hidden layers. A graphical representation of the flow of neural networks is depicted in Figure \ref{fig:deepnet_generic}.

    
For a neural network, we denote by $L$ the number of hidden layers (termed as depth) and by $N$, the maximum number of neurons at the hidden layers (denoted by width). Henceforth, we denote by $\cF_{NN}(d, N, L, W, o)$ by the collection of all neural networks with depth $L$, width $N$, the total number of active (non-zero) weights $W$, input dimension $d$ and output dimension $o$. We may sometimes omit the input and output dimensions in the specification of the class of neural networks when it is clear from the context. 

The expressibility of neural networks has been a topic of interest for decades, which started to flourish at the beginning of 1990s. It was proved in \cite{hornik1989multilayer} that the set of all shallow neural networks (i.e. with one hidden layer) is dense in the space of all Borel measurable functions and any Borel measurable function $f$ can be estimated within arbitrary accuracy by increasing the number of neurons in the hidden layer. The degree of accuracy was quantified by \cite{barron1993universal} -- if the Fourier transform of the function has a finite first moment, then the approximation error decreases at a rate $1/\sqrt N$, where $N$ is the number of hidden neurons. The error of approximating any \emph{smooth} function by a multi-layer feed-forward network was established by \cite{mhaskar1996neural} along with a choice of $N, L$ which achieves the optimal accuracy. Since then, there has been a series of research on the approximation capability of deep neural networks in terms of their width, height, weights, and activation function -- recently, researchers have obtained precisely bound on the approximation error of smooth functions in terms of depth and width of neural network ~cf. \cite{yarotsky2017error, lu2021deep, fan2022factor} and references therein. For example, Theorem 1.1 of \cite{lu2021deep} shows that if $f:[0, 1]^d \to \reals$ is $\beta$ times differentiable with bounded derivatives then there exists a neural network $\phi$ with width $c_1 N\log{N}$ and depth $c_2 L\log{L} + d$ ($c_1, c_2$ are some constants depending on $\beta$) with the following approximation error: 
$$
\|f - \phi\|_\infty \le C(NL)^{-\frac{2\beta}{d}} \,,
$$
where the constant $C$ depends on $(d, \beta)$. 

A related line of research delves into understanding how DNNs can successfully adapt to the unknown underlying low dimensional structure of the mean function in nonparametric regression. Performing nonparametric regression using suitable DNNs and achieving minimax optimal error bound (up to some polylog factor) initially started from the pioneering work of \cite{kohler2005adaptive}. The minimax rate of estimation of a mean function $f$ in a standard nonparametric regression problem (e.g., additive noise model $Y = f(X) + \eps$) is $n^{-2\beta/(2\beta +d)}$ where $\beta$ is the smoothness index of $f$ (Assumption \ref{assn:smoothness}) and $d$ is the dimension of $X$. The estimation error suffers from the curse of dimensionality, i.e. the rate is very slow for large $d$. However, as mentioned in the introduction, it is possible to circumvent this curse by imposing more structure on $f$ such as  additive (see \cite{stone1985additive}), or some higher order interaction model (see \cite{stone1994use}). 
In \cite{kohler2005adaptive}, the authors show that it is possible to adapt to the rate $n^{-2\beta/(2\beta + k)}$ using neural network-based estimate if the model is a $k$-way interaction model. The key advantage of using a neural network is that we only need to specify its architecture, i.e. the width and layer, not the exact structure of the mean function and the estimated NN becomes minimax rate optimal. Recently, in a series of papers \citep{kohler2016nonparametric, schmidt2020nonparametric, bauer2019deep, kohler2021rate, fan2022noise, fan2022factor}, the authors also establish that it is possible to estimate the underlying mean function $f_0$ at a rate $n^{-2\beta/(2\beta + k)}$ via neural networks when $f_0$ belongs to a more general class than $k$-way interaction model, called \emph{generalized hierarchical interaction model}. 

All the above analyses assume that the underlying dimension of $d$ is fixed. This is at odds with many modern applications.  Often time the constant in front of the rate of estimation depends on the dimension of $X$ even for the simple additive model as we need to estimate $d$ univariate functions. When $d$ is fixed, it is possible to get away with a cruder bound. However, if $d \rightarrow \infty$, we need new techniques to obtain the optimal dependence on $d$, which is  one of the primary goals of this paper.  

\begin{figure}
    \centering
    \includegraphics[scale=0.4]{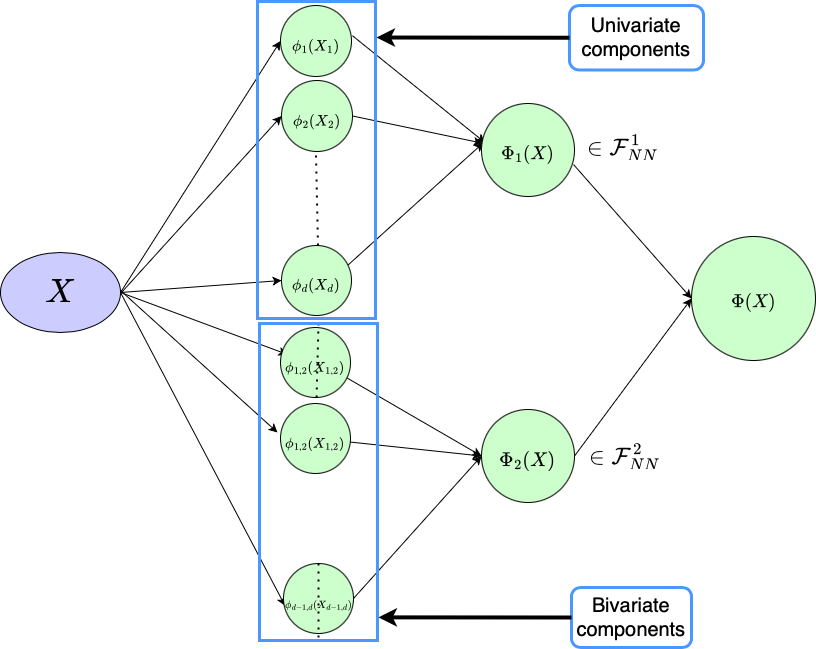}
    \caption{A schematic diagram for the structured deep neural networks that are used to estimate the structured nonparametric regression with interactions. The first part consists of additive fully connected neural networks for the univariate predictors and the second component comprises the summation of fully connected bivariate neural networks for approximating bivariate interaction effects.}
    \label{fig:deepnet_addition}
\end{figure}

 We now briefly describe our estimation procedure. Consider the model in equation \eqref{eq:two_way_ss_anova}. We select two neural networks $(\hat \phi_1, \hat \phi_2)$ via minimizing the squared-error loss: 
\begin{equation}\label{eq:estimator_growing_d}
\left(\hat \phi_1, \hat \phi_2\right) = \argmin_{\phi_1 \in \cF_{NN}^1,\phi_2 \in \cF_{NN}^2}\frac1n \sum_{i=1}^n \left(Y_i - \phi_1(X_i) - \phi_2(X_i)\right)^2
\end{equation}
{\color{black}and set estimate $\hat f_{\text{big}} := \hat \phi_1 + \hat \phi_2$. Finally set the estimator $\hat f : = sgn (\hat f_{\text{big}})(|\hat f_{\text{big}}| \wedge B)$, i.e. truncate the output of $\hat f$ at $[-B, B]$ as we assume $\|f_0\|_\infty \le B$ (see point 3 of Assumption \ref{assn:identifiability})}. Here the class of neural networks $\cF^1_{NN}$ and $\cF^2_{NN}$ are defined as: 
\begin{align}
    \label{eq:def_F1} \cF_{NN}^1 & = \cF_{NN}\left(d, c_1 dN_1\log{N_1}, c_2L_1\log{L_1}, c_3 d(N_1\log{N_1})^2L_1\log{L_1}, 1\right) \\
    \label{eq:def_F2} \cF_{NN}^2 & = \cF_{NN}\left(d, c_1 d^2N_2\log{N_2}, c_2L_2\log{L_2}, c_3 d^2(N_2\log{N_2})^2L_2\log{L_2}, 1\right)
\end{align}
for some constants $c_1, c_2, c_3$. For example, in $\cF^2_{NN}$, all pairwise interactions are used as input variables.  Note that we are not using a fully connected neural network here. Rather, we estimate each univariate (resp. bivariate) component function via fully-connected neural networks of width $N_1 \log{N_1}$ (resp. $N_2\log{N_2}$) and depth $L_1 \log{L_1}$ (resp. $L_2\log{L_2}$) and then add them.
See Figure~\ref{fig:deepnet_addition} for an illustration.  With slight abuse of notation, we still use $\cF_{NN}^1$ and  $\cF_{NN}^2$ to denote these specially structured neural networks.
The quantities $N_i$ and $L_i$ typically depend on the sample size $n$ and will be specified later. Like any other non-parametric estimator, the \emph{generalization error} of $\hat f$ can be decomposed into two parts: 
\begin{equation}\label{eq:error_decomp}
    \|\hat f - f_0\|_{L_2(P_X)} \le \underbrace{\|\hat f - f^\star\|_{L_2(P_X)}}_{\text{Statistical error}}  + \underbrace{\|f^\star - f_0\|_{L_2(P_X)}}_{\text{Approximation error}} 
\end{equation}
Here $f^\star = \phi_1^\star + \phi_2^\star$ is (approximately) the best approximator of $f_0$ in terms of $L_2(P_X)$ norm among the class of neural networks, i.e. 
$$
\|f^\star - f_0\|_{L_2(P_X)} \le \inf_{\phi_1 \in \cF_{NN}^1, \phi_2 \in \cF_{NN}^2}\|f_0 - \phi_1 - \phi_2\|_{L_2(P_X)} + \tau_n
$$
for some small tolerance $\tau_n$ to be specified later. 
The rest of our analysis is divided into two parts: Subsection \ref{sec:approx_lowdim} bounds the approximation error and Subsection \ref{sec:stat_lowdim} bounds the statistical error in terms of $(d, N_i, L_i, n)$. Finally, in Subsection \ref{sec:minimax_lowdim} we choose the architectures of $\cF_{NN}^1, \cF_{NN}^2$ to balance the errors. This will prove that the rate obtained by the neural network is minimax optimal both in terms of $(n, d)$, up to some logarithmic factor. 
\begin{remark}
\label{rem:diff_arch}
    It is important to note, we are estimating $f_0$ as the sum of two neural networks with different architectures instead of using a single neural network. To see why this is necessary, note that the rate of convergence of the estimator is $dn^{-2\beta_1/(2\beta_1 + 1)}$ and $d^2n^{-\beta_2/(\beta_2 + 1)}$ (Corollary \ref{cor:error_lowdim}) upto a logarithmic factor, where $\beta_1$ is the smoothness index of the univariate component and $\beta_2$ is the smoothness index of the bivariate components of $f_0$. Now note that if $\beta_2$ is sufficiently large, then $dn^{-2\beta_1/(2\beta_1 + 1)}$ will be dominant rate, otherwise $d^2n^{-\beta_2/(\beta_2 + 1)}$ will be the dominant one. If we only use one neural network, then the rate obtained is {\color{black}$d^2(n^{-2\beta_1/(2\beta_1 + 1)}  + n^{-\beta_2/(\beta_2 + 1)})$}, which is not minimax optimal (Theorem \ref{thm:minimax_lower_bound_growing_d}) especially when $\beta_2 \gg \beta_1$. The trade-off between the error terms is more complicated as $d:=d(n) \rightarrow \infty$ and using two neural networks (one for univariate components and the other for bivariate components) we can obtain the optimal error bounds.
    We conjecture that it is not possible to obtain the minimax error bound by using a single neural network. 
\end{remark}

\subsection{Approximation error}
\label{sec:approx_lowdim}
In this section, we compute the approximation error of mean estimation in $2$-way interaction model \eqref{eq:two_way_ss_anova}.
To this end, we invoke Theorem 1.1 of \cite{lu2021deep}, which states any $\beta$ smooth function $f:[0, 1]^d \to \reals$ can approximate $f$ within error of $(NL)^{-2\beta/d}$ in $L_\infty$ norm using a DNN with width $O(N\log{N})$ and depth $O(L\log{L} + d)$. Recall that, by Assumption, \ref{assn:identifiability}, $\int f_j(x)dx=0$ and both the marginal integrals of $f_{ij}$s are $0$. As stated in the previous section, we will approximate the $f_i$s and $f_{ij}$s by two separate neural networks with different depths and widths. Before stating our main theorem, we first quantify a bound on the growth of $d$ in comparison to $n$ in the following assumption: 
\begin{assumption}
    \label{asm:d_growth}
    We assume that $d$ is growing with $n$ and 
    $$
    \left(dn^{-\frac{2\beta_1}{2\beta_1 + 1}} + \binom{d}{2}n^{-\frac{\beta_2}{\beta_2 + 1}}\right)\log^{4.5}{n} = o(1) \,.
    $$
\end{assumption}
The following theorem establishes the approximation error of the mean function via neural networks: 
\begin{theorem}
\label{thm:approx_growing_d}
    Consider the two-way interaction model in the equation 
    \eqref{eq:two_way_ss_anova}, 
    where all the univariate components are $\beta_1$ smooth and all the bivariate components are $\beta_2$ smooth. {\color{black} Choose $N_1 L_1 = \lfloor n^{1/2(2\beta_1 + 1)}\rfloor$, $N_2 L_2 = \lfloor n^{1/2(\beta_2 + 1)}\rfloor$} Then we have: 
    \begin{align*}
    & \inf_{\phi_2 \in \cF^1_{NN}, \phi_2 \in \cF^2_{NN}}\bbE\left[\|f_0(X) - \phi_1(X) - \phi_2(X)\|^2\right] \\
    & \qquad \qquad \qquad \le C_3 \left(d(N_1L_1)^{-4\beta_1} + \binom{d}{2} (N_2L_2)^{-2\beta_2}\right) \,,
    \end{align*}
   where $\cF^1_{NN}$ and $\cF^2_{NN}$ are same as defined in \eqref{eq:def_F1} and \eqref{eq:def_F2}.
   All the constants $(C_1, C_2, C_3)$ are independent of $(d, N_1, L_1, N_2, L_2)$. 
\end{theorem}

\noindent The proof of Theorem \ref{thm:approx_growing_d} is deferred to Section \ref{sec:proof_approx_growing_d}. Here we provide the sketch of the proof. The key idea in our proof is to construct a neural network that not only approximates $f$ but also (approximately) satisfies the identifiability constraint Assumption \ref{assn:identifiability}, i.e. the marginals of $f_{ij}$'s are $0$. This will nullify the effect of the higher-order terms yielding the optimal rate of approximation. For univariate components, this integral constraint can be satisfied exactly; consider a function $f:[0, 1] \to \reals$ where $f$ satisfies Assumption \ref{assn:smoothness} with some $\beta$ and $\int_0^1 f(x) \ dx = 0$. Then by \cite[Theorem 1.1]{lu2021deep}, there exists a neural network $\phi$ with width of $O(N_1\log{N_1})$ and depth $O(L_1\log{L_1})$ such that $$\|f - \phi\|_\infty \le C(N_1L_1)^{-2\beta}$$ where the constant $C$ depends only on $\beta$ and Sobolev norm of $f$. 
Define $I_\phi := \int_0^1 \phi(x) \ dx$ and define a new neural network $\tilde \phi = \phi - I_\phi$ guaranteeing $\int \tilde \phi=0$. As subtracting a constant from a neural network amounts to changing the bias of the last layer, the architecture of $\phi$ and $\tilde \phi$ are the same. Using triangle inequality,
$$
\|f - \tilde \phi\|_\infty \le 2C(N_1L_1)^{-2\beta} \,.$$
The above idea doesn't work for the approximation of the bivariate functions $f_{ij}$ as we not only need the integral to be $0$ but also the marginals to be $0$, i.e. $\int f_{ij}(x, y) \ dx = \int f_{ij}(x, y) \ dy = 0$. To address this issue, we devise a debiasing technique summarized as follows: given any bivariate function $f$ which is $\beta$ times differentiable with bounded derivatives, there exists a neural network $\phi(x, y)$ (using \cite[Theorem 1.1]{lu2021deep}) with architecture $O(N_2\log{N_2}), O(L_2\log{L_2})$ that satisfies: 
$$
\|f - \phi\|_\infty \le C(N_2L_2)^{-\beta} \,.
$$
Given such a neural network $\phi$, we construct $\phi_1$ (resp. $\phi_2$) such that $\phi_1(x)$ (resp. $\phi_2(y)$) is a neural network with width $O(N_2\log{N_2})$ and depth $O(L_2\log{L_2})$ 
which approximates $\int \phi(x, y) \ dy$ (resp. $\int \phi(x, y) \ dx$). Our final estimator is of the form $\tilde \phi = \phi - \phi_1 - \phi_2 + \iint \phi(x, y) \ dx \ dy$. We show that $\|\tilde \phi- f\|_{\infty}=O((NL)^{-\beta})$ and $\int \tilde \phi(x) \ dx = \int \tilde \phi(x, y) \ dy = O\left((NL)^{-2\beta}\right)$. 
This smaller approximation error ($(N_2L_2)^{-2\beta}$ instead of $(N_2L_2)^{-\beta}$) is crucial to obtain the optimal error bounds. Details of the proof can be found in Appendix \ref{sec:proof_approx_growing_d}. 

\subsection{Statistical error and main theorem}
\label{sec:stat_lowdim}
In this section, we discuss the techniques to bound the statistical error of estimation and also present the main theorem that bounds the out-of-sample prediction error. Statistical error of an estimator typically relies on the number of samples and the complexity of the underlying class of function. Recall from \eqref{eq:estimator_growing_d} that we approximate $f_0$ by the sum of two neural networks $\hat \phi_1 + \hat \phi_2$ where $\hat \phi_1$ approximates the univariate component and $\hat \phi_2$ approximates the bivariate component. It is immediate from the proof of Theorem \ref{thm:approx_growing_d} that, the approximation error depends on $N$ and $L$ through their product $NL$. Therefore, we work with the neural network with constant depth and optimize over width. So, $\cF^1_{NN}$ and $\cF^2_{NN}$ (equation \eqref{eq:def_F1} and \eqref{eq:def_F2}) can be revised as: 
\begin{align}
    \label{eq:modified_def_F1} \cF_{NN}^1 & = \cF_{NN}\left(d, c_1 dN_1\log{N_1}, c_2, c_3 d(N_1\log{N_1})^2, 1\right), \\
    \label{eq:modified_def_F2} \cF_{NN}^2 & = \cF_{NN}\left(d, c_1 \binom{d}{2}N_2\log{N_2}, c_2, c_3 \binom{d}{2}(N_2\log{N_2})^2, 1\right).
\end{align}
In practice, one estimates $\hat \phi_1, \hat \phi_2$ by optimizing over the weights of the neural network through gradient descent (e.g. see 6.2 of \cite{goodfellow2016deep}). Since our problem is highly non-convex in terms of the weights of the neural network, it is not immediately clear whether such an algorithm will converge to global minima. To avoid such issues, we work with a neural network that minimizes empirical risk in our article. Understanding the properties of the gradient descent-based estimator is left open for future research. 

The aim of this section is to bound the \emph{variance} term, i.e. $\bbE[(\hat f(X) - \phi^\star(X))^2 \mid \bS_n]$ where $\phi^\star$ is the best approximator of the mean function $f_0$ and $\bS_n:= \{(X_i,Y_i), i \in [n]\}$ be the training sample. The prediction/test error depends on the trade-off between the complexity of the underlying function class and the number of training samples. If the underlying function class is too complex compared to the number of samples $n$, then it may overfit the data resulting in high generalization error. Therefore, it is imperative to quantify the complexity of the class of neural networks. The following lemma of \cite{bartlett2019nearly} establishes a bound on the VC dimension of the class of neural networks with the number of active weights $W$ and depth $L$:  
\begin{lemma}[(5), \cite{bartlett2019nearly}]\label{lem:bartlett_vc_bound}
    Suppose $V(W,L)$ denotes the largest VC- dimension of a ReLU network with $W$ parameters and $L$ layers. There exist constants $c_1, c_2>0$ such that
    \[ c_1 WL \log(WL) \le V(W,L) \le c_2 WL \log{WL}. \]
\end{lemma}
\noindent Roughly speaking, if we have a fully connected neural network with width $N$ and depth $L$, then $W \sim N^2L$ and consequently the VC dimension is $\sim N^2L^2$.  
Once we have a bound on the complexity of the underlying function class, we can use techniques from empirical process theory to obtain the rate of convergence of the generalization error. For a VC-type function class with VC dimension $V$, the statistical error is $O(\sqrt{V/n})$ up to logarithmic factors. 
\noindent To see the necessity of the assumption, consider the rate obtained in Corollary \ref{cor:error_lowdim}. Assumption \ref{asm:d_growth} ensures the consistency of the estimator. On the other hand, it is also necessary up to the logarithmic factor as we establish in Section 4 that the rate obtained in Corollary \ref{cor:error_lowdim} is minimax optimal up to the logarithmic factor (Theorem \ref{thm:minimax_lower_bound_growing_d}). Therefore the assumption is the weakest possible assumption on the growth of $d$ with respect to $n$ under no further structural assumption. The following theorem presents a bound on the overall approximation error, combining both statistical and approximation errors.
\begin{theorem}[Main theorem]
    \label{thm:statistical_error_lowdim}
    Under Assumption \ref{assn:density_bound} - \ref{assn:smoothness} and \ref{asm:d_growth}, the ERM estimator $\hat f$ of $f_0$ defined in equation \eqref{eq:estimator_growing_d} satisfies: 
    $$
    \bbE\left[\left(\hat f(X) - f_0(X)\right)^2 \mid \bS_n\right] = O_p\left(\rho_n^2 + \frac{V_n}{n}\log^{3/2}{n}\right),
    $$
    where
    \begin{align*}
        \rho_n^2 &= \text{approximation error} \le C_1\left(dN_1^{-4\beta_1} + \binom{d}{2}N_2^{-2\beta_2}\right) \\
        V_n &= \text{complexity} \le C_2\left(dN_1^2\log^2{N_1}\log{(dN_1)} + \binom{d}{2}N_2^2\log^2{N_2}\log{(dN_2)}\right) \,.
    \end{align*}
\end{theorem}
\noindent
The proof of the theorem can be found in Appendix \ref{sec:proof_statistical_error_lowdim}. Now we need to choose the value $N_1, N_2$ carefully to balance both the approximation error and statistical error. This leads to the following corollary: 
\begin{corollary}
\label{cor:error_lowdim}
Choosing $N_1 = n^{1/2(2\beta_1+ 1)}$ and $N_2 =n^{1/2(\beta_2 + 1)}$ (take the nearest integer if they are not integers) in Theorem \ref{thm:statistical_error_lowdim}, we have: 
$$
\bbE\left[\left(\hat f(X) - f_0(X)\right)^2 \mid \bS_n\right] =O_p\left(\left(dn^{-\frac{2\beta_1}{2\beta_1 + 1}} + \binom{d}{2}n^{-\frac{\beta_2}{\beta_2 + 1}}\right)\log^{4.5}{n}\right) \,.
$$ 
\end{corollary}
\noindent In the subsequent section, we show that this rate is minimax optimal up to the logarithmic factor. An immediate question that stems from this analysis is whether we can improve the dependence on the logarithmic factor. The recent article \citep{fan2022factor} proposes a slightly different approximation technique than that of \cite{lu2021deep}, which might reduce the power of the logarithmic factor from $4.5$. 
However, the more interesting question is whether we can completely get rid of the log factor which is still open. 
\begin{remark}\label{rem:extension_to_k}[Extension for k-way interaction model]
Our analysis for the two-way interaction model can also be extended verbatim to the $k$-way interaction model. Here we provide a sketch of the extension. For identifiability of the model, we need to assume all marginals of $j$-component functions are $0$, $j \ge 2$. We estimate $f_0$ via the sum of $k$ neural networks $\phi_1 + \dots + \phi_k$ where 
$$
\phi_i \in \cF_{NN}\left(d, c_1\binom{d}{i}N_i\log{N_i}, c_2, \binom{d}{i}(N_i\log{N_i})^2, 1\right) \,.
$$
Since the debiasing technique used in our proofs rely on the approximation of polynomials by neural networks, they can be extended to a general $k$-way interaction model. This in turn controls the approximation error. The statistical error can be bounded similarly by the VC-dimension of the neural networks via Lemma \ref{lem:bartlett_vc_bound}. 
\end{remark}

\subsection{Minimax lower bound}
\label{sec:minimax_lowdim}
In this section, we establish the minimax lower bound for estimating the non-parametric $2$-way interaction model \eqref{eq:two_way_ss_anova}, when the dimension of the covariate $d= d(n) \rightarrow \infty, d=o(n)$ and satisfies Assumption \ref{asm:d_growth}. 
For our estimation problem, the minimax risk is defined as: 
$$
\mathfrak{M}(n, d, \cF) = \inf_{\hat f} \sup_{\substack{f \in \cF \\ X \sim P_{X}}}\bbE_f\left[\left(\hat f(X) - f(X)\right)^2\right]
$$
Here $f(X) = \sum_{j=1}^d f_j(X_j) + \sum_{i < j} f_{ij}(X_i, X_j)$ and the supremum is taken over all $\{f_j\}$ and $\{f_{ij}\}$ where the component functions belongs $\Sigma(\beta, L)$ (see Assumption \ref{assn:smoothness}). 
\\\\
\noindent
The main difficulty in establishing the minimax lower bound is to incorporate the effect of growing dimension $d$. Our problem is certainly harder than the problem of estimating $f_1(x_1) + f_{12}(x_1, x_2)$, assuming all other components are known in advance.  The issue is how other components contribute to the lower bound in the growing $d$ setting.  The following theorem gives a precise answer.
\\\\
\begin{theorem}
\label{thm:minimax_lower_bound_growing_d}
    Consider the two-way interaction model as defined in equation \eqref{eq:two_way_ss_anova} where each component functions $\{f_i\}_{i=1}^d$ and $\{f_{ij}\}_{1 \le i < j \le d}$ belongs to $\Sigma(\beta, L)$ (see Assumption \ref{assn:smoothness}) and denote this collection of mean functions as $\cF$. Then the minimax rate of estimation under Assumptions \ref{assn:density_bound}-\ref{assn:smoothness} is: 
    $$
    \mathfrak{M}(n, d, \cF) = \inf_{\hat f} \sup_{\substack{f \in \cF \\ X \sim P_{X}}}\bbE_f\left[\left(\hat f(X) - f(X)\right)^2\right] \ge c\left( dn^{-\frac{2\beta_1}{2\beta_1 + 1}} +  \dbinom{d}{2}n^{-\frac{2\beta_2}{2\beta_2 + 2}}\right) \ \,.
    $$
    for some constant $C$ independent of $(n, d)$.  
\end{theorem}
\noindent
The proof of this theorem can be found in Appendix \ref{sec:proof_minimax_growing}. A few remarks are in order; first of all, Theorem \ref{thm:minimax_lower_bound_growing_d} establishes the fact that neural network-based estimate of the mean function $f_0$ is minimax rate optimal up to a poly-log factor. This result complements the one in \cite{kohler2021rate, schmidt2020nonparametric} for the non-parametric regression in a fixed dimension regime. Secondly, Theorem \ref{thm:minimax_lower_bound_growing_d} is derived under Assumption \ref{assn:smoothness}, where we assume all univariate components are $\beta_1$-smooth and all bivariate components are $\beta_2$-smooth. 
However, our proof can be easily adapted to the setting where each function has different smoothness and we have to pay the price for the lowest smoothness. See Remark \ref{rem:extension_to_k} where we have pointed out the key steps to prove such extension.  

\begin{remark}\label{rmk:general_lbd}[Extension for k-way interaction model]
The proof of the minimax lower bound for the two-way interaction model can be easily generalized for the general $k-$way interaction model by constructing alternatives similarly and using Fano's inequality. The optimal rate is given by 
$$
\inf_{\hat f} \sup_{f \in \cF, X \sim P_{X}} \bbE_{P_{X}}[(\hat f(X) - f(X))^2] \gtrsim \left(\sum_{j=0}^k \dbinom{d}{j}n^{-\frac{2\beta_j}{2\beta_j + j}}\right)
$$
\end{remark}

\begin{remark}
We have studied the effect of dimensionality in the minimax upper and lower bound for the interaction model in our article. A natural next step is to find Pinsker's constant for the $k$-way interaction model in diverging dimensions, which is left as a future research direction.
\end{remark}
\section{Analysis of DNN when $d \gg n$}\label{sec:high_dim_results}
This section presents our analysis when the underlying dimension of $X$ is larger than the sample size. Although simpler parametric models have been well-studied to avoid the curse of dimensionality, the literature on nonparametric estimation for $k$-way interaction model \eqref{eq:model_k} in this regime is relatively sparse. In the nonparametric framework, a significant amount of research has been done on a high dimensional additive model under sparsity constraint where the model under consideration is as follows: 
$$
Y_i = \sum_{j \in S}f_{0, j}(X_{ij}) + \eps_i
$$
where $S \subseteq \{1, \dots, d\}$ with $s = |S| \ll n$. This is a standard assumption in high dimensional statistical analysis, where the true signal depends on a few covariates, but the active set $S$ is apriori unknown. The analysis of the sparse additive model in high dimension was initiated in \cite{lin2006component} and later refined in a series of papers literature~\citep{koltchinskii2010sparsity,ravikumar2009sparse,tan2019doubly,yuan2016minimax}. \cite{koltchinskii2010sparsity} obtained error bounds under a global boundedness assumption, which was subsequently removed by \cite{raskutti2012minimax}. Recently, \cite{tan2019doubly} obtained minimax guarantees when the underlying component functions lie in a reproducing kernel Hilbert space. All the previous works typically use two penalties for optimal estimation: one to control the complexity of the underlying function class and the other to control the sparsity. Yet, it is enough to use only one penalty for inducing sparsity for a neural network-based estimator, as the complexity can be controlled through the network's width and depth. 

In spite of such extensive works on sparse linear models, the study of the high-dimensional linear models with $k$-way interaction ($k>1$) is very sparse, letting alone $k$-way non-parametric interactions.   For some references about the variable selection in high dimensional linear interaction models, readers may consult \cite{zhao2009composite}, \cite{bien2013lasso}, \cite{hao2014interaction} and references therein. In this section, we bridge this gap by studying the asymptotic properties of the neural network estimator.
\par As before, we only elaborate on the analysis of the $2$-way interaction model and comment on how to extend our analysis for the general $k$-way interaction model. The two-way sparse additive model is defined as follows: 
\begin{equation}\label{eq:define_sparse_model}
    Y_i= \sum_{j \in S_1}f_{0, j}(X_{ij})+ \sum_{(k,l)\in S_2} f_{0, kl}(X_{ik},X_{il}) + \varepsilon_i \equiv f_0(X_i)+\varepsilon_i,
\end{equation}
where $S_1 \subset [d]$, $S_2 \subset [d]\times[d]$ with $S_1,S_2$ sparse, i.e., $s_i:=|S_i| \ll n$ for $i=1,2$. Define this class of functions by $\mathcal{F}_{\text{sp}}$. Here also we assume that the univariate components are $\beta_1$ smooth and the bivariate components are $\beta_2$ smooth. From Theorem \ref{thm:statistical_error_lowdim}, there exists $\{\phi^\star_j\}_{j \in S_1}, \left\{\phi^\star_{kl}\right\}_{(j, k) \in S_2}$, $j\in S_1, (j,k) \in S_2$, such that, $\phi^\star_j \in \mathcal{F}_{NN,1}$, $\phi^\star_{jk} \in \mathcal{F}_{NN,2}$ and
\begin{align*}
    & \|f_{0,j} - \phi^\star_{j}\|_\infty \le C_1 N_1^{-2\beta_1} \ \forall \ j \in S_1 \,,  \\
    & \|f_{0,kl} - \phi^\star_{kl}\|_\infty \le C_2 N_2^{-\beta_2} \ \forall \ (k < l) \in S_2 \,.
\end{align*}
where 
\begin{align}
    \label{eq:def_nn_uni}\cF_{NN, 1} & = \cF_{NN}(1, c_1N_1\log{N_1}, c_2, c_3 N^2_1\log^2{N_1}, 1) \\
    \label{eq:def_nn_bi} \cF_{NN, 2} & = \cF_{NN}(2, c_4 N_2\log{N_2}, c_5, c_6 N^2_2\log^2{N_2}, 1) \,.
\end{align}
Finally, define 
\begin{align}\label{eq:nn_approx_sparse}
    \phi^\star=\sum_{j \in S_1}\phi^\star_j +\sum_{(k<l) \in S_2} \phi^\star_{kl}.
\end{align}

Following Assumption \ref{asm:d_growth}, we assume the following growth condition on the sparsity parameter $s_1, s_2$: 
\begin{assumption}
    \label{asm:d_growth_high}
    The sparsity parameters $s_1, s_2$ corresponding to the univariate components and bivariate components satisfy the following condition: 
    $$
    s_1\left(n^{-\frac{2\beta_1}{2\beta_1 + 1}}\log^4{n} + \frac{\log{d}}{n}\right) + s_2\left(n^{-\frac{\beta_2}{\beta_2 + 1}}\log^4{n} + \frac{\log{d}}{n}\right) = o(1) \,.
    $$
\end{assumption}
The reason behind Assumption \ref{asm:d_growth_high} is similar for that of Assumption \ref{asm:d_growth}, i.e. without this assumption, there will be no consistent estimator (modulo the logarithmic factor) as it is minimax optimal rate of estimation. The analysis for the high dimensional model is not a straightforward extension of the techniques used when $d = o(n)$. To understand why, recall that even in a simple sparse linear model with $\ell_1$ penalty, it is not possible to obtain a minimax rate optimal estimator without a form of \emph{restricted strong convexity} assumption, which is typically not needed when the dimension grows slower than the sample size. Therefore, for ease of presentation, we divide the entire analysis into three parts: Section \ref{sec:fixed_design_highdim} establishes the rate of convergence for the fixed design model, i.e. when $X_i$'s are some fixed points in $\reals^p$. Section \ref{sec:random_design_highdim} deals with the random design, i.e. $X$ is assumed to be a random variable. Finally, in Section \ref{sec:minimax_high_dim} we present the minimax lower bound to establish that our neural network-based estimator is minimax rate optimal up to log factors.

\subsection{Analysis of Fixed design model}
\label{sec:fixed_design_highdim}
In this section, we present our analysis for the fixed design model, where we assume $X_i$'s to be fixed. Similar to \eqref{eq:estimator_growing_d}, we estimate univariate and bivariate components using neural networks with different architecture, but at the same time, we use the $\ell_1$-norm of $L_2(\bbP_n)$ penalty to enforce sparsity.
See page 966 of \cite{antoniadis2001regularization} and \cite{yuan2006model} for such an idea in selecting a group of variables. 
Our estimator $\hat f = \sum_j \hat \phi_j + \sum_{k < l} \hat \phi_{kl}$ where components are defined as 
\begin{align}
\label{eq:fixed_sparse_penalize}
\left\{\hat \phi_j\right\}, \left\{\hat \phi_{kl}\right\} & = \mathop{\arg \min}\limits_{\substack{\phi_j \in \cF_{NN, 1} \\ \phi_{kl} \in \cF_{NN, 2}}} \left[\frac1n \sum_i \left(Y_i - \sum_j \phi_j(X_{ij}) - \sum_{k < l}\phi_{kl}(X_{ik}, X_{il})\right)^2 \right.\nonumber \\
& \qquad \qquad \qquad \qquad \qquad \left. + \left(\lambda_{n,1}\sum_j \|\phi_j\|_n + \lambda_{n,2}\sum_{k < l} \|\phi_{kl}\|_n\right)\right]
\end{align}
where $\cF_{NN, 1}$ and $\cF_{NN, 2}$ are defined in \eqref{eq:def_nn_uni} and \eqref{eq:def_nn_bi} respectively, $\| \cdot \|_n$ is the $L_2(\bbP_n)$ norm and $\lambda_{n,1}, \lambda_{n,2}$ to be specified later.\par
We now make couple of remarks on the estimation procedure. First, although we put $L_2(\bbP_n)$ penalty on the component functions to enforce sparsity, we may use the $\ell_1$ penalty on the last layer of the weights (which adds the component neural networks). Such a penalty, in spite of its computational benefit, is difficult to analyze theoretically. Second, as we do not necessarily assume $\beta_1=\beta_2$, we need two regularization parameters $\lambda_{n, 1}$ and $ \lambda_{n, 2}$. We also need different architectures to adapt to the smoothness as highlighted by Remark \ref{rem:diff_arch}. Furthermore, for technical simplicity, we assume the following $l^\infty$ bound: 
\begin{equation}\label{eq:B_bound}
    \|f_j\|_{\infty} \le B, \qquad \|f_{jk}\|_{\infty}\le B.
\end{equation}
Consequently, for numerical stability, we truncate the component neural networks (i.e. $\phi_j, \phi_{jk}$) at level $B$.
The following theorem establishes the rate of convergence.
\begin{theorem}
    \label{thm:fixed_design_wo_rsc}
    Under Assumption \ref{assn:density_bound}-\ref{assn:smoothness} and \ref{asm:d_growth_high} along with \eqref{eq:B_bound},
    the estimator obtained in \eqref{eq:fixed_sparse_penalize}  satisfes
    $$
    \|\hat f - f_0\|_n^2 = O_p\left(s_1\left(\rho_{n, 1}^2 +\lambda_{n, 1} + \frac{\lambda_{n,1}^2}{2}\right)  + s_2 \left(\rho_{n, 2}^2  + \lambda_{n,2} + \frac{\lambda_{n, 2}^2}{2}\right)\right)
    $$
    where $\rho_{n, 1}$ (resp. $\rho_{n, 2}$) is the approximation error of the univariate (resp. bivariate) components of the mean function by neural networks, bounded by \eqref{eq:approx-err}, provided that 
    \begin{equation}\label{eq:define_lambda}
    \lambda_{n, 1} = C_3\sqrt{\frac{V_{n, 1} \log{n}}{n} + \frac{2\log{d}}{n}}, \quad 
\lambda_{n, 2} = C_4\sqrt{\frac{V_{n, 2} \log{n}}{n} + \frac{3\log{d}}{n}}   \,,
\end{equation}
with $V_{n, 1}  = N_1^2 \log^3{N_1}$ and $V_{n, 2} = N_2^2 \log^3{N_2}$.
\end{theorem}
\noindent 
Let us try to understand and simplify the rate obtained by Theorem \ref{thm:fixed_design_wo_rsc}. Recall from Theorem \ref{thm:approx_growing_d}, we have:  
\begin{equation} \label{eq:approx-err}
\rho_{n, 1} \le C_1 N_1^{-2\beta_1}, \quad \rho_{n, 2} \le C_2 N_2^{-\beta_2} \,,
\end{equation}
for some constants $C_1, C_2>0$.  
Furthermore, it is revealed in our proof (see  \eqref{eq:lambda_1_val} and \eqref{eq:lambda_2_val}) that  $\lambda_{n, 1}$ and $\lambda_{n, 2}$ given by \eqref{eq:define_lambda} are the optimal choices, in which $V_{n, 1}$ and
$V_{n, 2}$ are the order of VC-dimensions for $\cF_{NN, 1}$ and $\cF_{NN, 2}$ respectively.
Now, typically $\lambda_{n, i} = o(1)$ (which holds as soon as $V_{n,i} = o(n/\log n)$, a condition necessary for consistency) and consequently $\lambda_{n,i}^2=o(\lambda_{n, i})$. Hence the rate of convergence in Theorem \ref{thm:fixed_design_wo_rsc} can be simplified as
$$
 \|\hat f - f_0\|_n^2 = O_p\left(s_1\left(\rho_{n, 1}^2 +\lambda_{n, 1} \right)  + s_2 \left(\rho_{n, 2}^2  + \lambda_{n,2}\right)\right) \,.
$$
Now we choose $N_1, N_2$ to balance $\rho_{n, i}^2$ and $\lambda_{n, i}$, which leads to the following corollary: 
\begin{corollary}
\label{cor:high_dim_rate}
Choosing $N_1 = n^{1/(2(1 + 4\beta_1))}$ and $N_2 = n^{1/(2(1 + 2\beta_2))}$ (take the nearest integer if they are not integers), we obtain: 
$$
 \|\hat f - f_0\|_n^2 = O_p\left(s_1\left(n^{-\frac{2\beta_1}{1 + 4\beta_1}}\log^2{n} + \sqrt{\frac{\log{d}}{n}}\right)  + s_2 \left(n^{-\frac{\beta_2}{1 + 2\beta_2}}\log^2{n} + \sqrt{\frac{\log{d}}{n}}\right)\right) \,.
$$
\end{corollary}
cAs we will see in Subsection \ref{sec:minimax_high_dim}, this rate is not minimax optimal. The reason is similar to that for the high dimensional sparse regression model: if we do not assume any condition on the curvature of the loss function (e.g. restricted isometry or restricted eigenvalue type conditions) then it is not possible to obtain any minimax optimal estimator which is computable in polynomial time as proved in \cite{zhang2014lower}. 
Hence, we need to assume a form of \emph{Restricted Strong Convexity} (RSC) to obtain a minimax optimal error bound. RSC-type assumptions were popularized by \cite{bickel2010hierarchical,candes2007dantzig}. A similar assumption is also used in the analysis of the sparse additive regression model, e.g. see \cite[Assumption 3]{tan2019doubly}. In particular, we assume the following:
\begin{assumption}
\label{assn:rsc}
There exist constants $\kappa_1, \kappa_2 > 0$, such that for any function $\phi = \sum_j \phi_j +\sum_{k<l} \phi_{kl}$ 
that satisfies:
\begin{align}
   & \lambda_{n, 1} \sum_{j \in S_1^c}\|\phi_j - \phi^\star_j \|_n + \lambda_{n, 2}\sum_{(k < l) \in S^c_2}\|\phi_{kl} - \phi^\star_{kl}\|_{n} \notag \\
    & \le 4 (s_1\rho_{n, 1}^2 + s_2 \rho_{n, 2}^2)  + 3\lambda_{n, 1} \sum_{j \in S_1}\|\phi_j - \phi^\star_j \|_{n} +  3\lambda_{n, 2} \sum_{(k < l) \in S_2}\|\phi_{kl} - \phi^\star_{kl}\|_{n} \notag \\
    \label{eq:cond_for_rsc} & \hspace{5em} + s_1\lambda_{n, 1}^2 + s_2\lambda_{n, 2}^2
 \end{align}
also satisfies: 
\begin{equation}
\label{eq:main_rsc_condition}
\kappa_1^2\sum_{j \in S_1}\|\phi_j - \phi^\star_j\|^2_{n} + \kappa_2^2 \sum_{(k < l) \in S_2}\|\phi_{kl} - \phi^\star_{kl}\|^2_{n} \le \|\phi - \phi^\star\|_n^2 
\end{equation}
where $\phi^\star$ is defined by \eqref{eq:nn_approx_sparse}, $\rho_{n, 1}$ and $\rho_{n,2}$ are defined in Theorem \ref{thm:fixed_design_wo_rsc} and $\lambda_{n, 1}$ and $\lambda_{n,2}$ are defined in \eqref{eq:define_lambda}.
\end{assumption}
\noindent 
Before going into further details, we first compare our assumption with the standard restricted eigenvalue (RE) condition for the high dimensional linear model \cite{bickel2009simultaneous}. Roughly speaking, the RE condition assumes: 
\begin{equation}
    \label{eq:re_explanation}
    \sum_{j \in S} (\beta_j - \beta_{0, j})^2 \lesssim \frac1n \|X(\beta - \beta_0)\|_2^2 
\end{equation}
for any $\beta$ that satisfies:  
\begin{equation}
    \label{eq:cone_cond}
    \sum_{j \in S^c} |\beta_j| \lesssim \sum_{j \in S}|\beta_j - \beta_{0, j}| \,,
\end{equation}
where $S$ is the true active set. Now consider the nonparametric additive model, where $X_j$ affects $Y$ through $f_j(X_j)$. Furthermore, we typically use the sum of $L_2(\bbP_n)$ norm as a penalty function in place of $L_1$ norm for sparsity. Ignoring the approximation error for the time being (i.e. assuming true $f_j$ belong to the function class over which we optimize), a natural analog of condition \eqref{eq:re_explanation} is the following: 
\begin{equation}
    \label{eq:re_explanation_nonparam}
     \sum_{j \in S} \|f_j - f_{0, j}\|_n^2  
     \lesssim \|\hat f - f_0\|_n^2 
\end{equation}
for all $f = \sum_j f_j$ that satisfies (similar to \eqref{eq:cone_cond})
\begin{equation}
\label{eq:cone_nonparam}
\sum_{j \in S^c}\|f_j \|_n \lesssim \sum_{j \in S} \|f_j - f_{0, j}\|_n \,.
\end{equation}
This is precisely what was used in \cite[Assumption 3]{tan2019doubly} for analyzing additive models in high dimensions. Assumption \ref{assn:rsc} is almost the same as this one with some modification due to the fact we are approximating the component functions through neural networks and consequently we may have non-zero approximation error. Ignoring the bivariate components, \eqref{eq:main_rsc_condition} is exactly the same as \eqref{eq:re_explanation_nonparam} with $f$ replaced by $\phi$ and $f_0$ replaced by $\phi^\star$, which is the best approximator of $f_0$ from the class of neural networks. The only difference is in 
\eqref{eq:cond_for_rsc}, which differs from \eqref{eq:cone_nonparam} as we need to take into account the approximation error. This is precisely why the term $s_1 \rho_{n, 1}^2 + s_2 \rho_{n, 2}^2 $ appears on the right-hand side of equation \eqref{eq:cond_for_rsc}. The other additional term $s_1 \lambda_{n, 1}^2 + s_2 \lambda_{n, 2}^2$ is due to some mathematical artifact of the proof and also can be made arbitrarily smaller in order. For simplicity, if we ignore the bivariate components and this additional term for the time being, then \eqref{eq:cond_for_rsc} simplifies to: 
\begin{align}
\label{eq:cond_for_rsc_additive}
   \lambda_{n, 1} \sum_{j \in S_1^c}\|\phi_j - \phi^\star_j \|_n  & \le 4 s_1\rho_{n, 1}^2 + 3\lambda_{n, 1} \sum_{j \in S_1}\|\phi_j - \phi^\star_j \|_{n} \,.
\end{align}
If the function class is well-specified (as in the case of linear regression or additive models in \cite{tan2019doubly}), then the first term of the RHS can be removed. Consequently, cancelling $\lambda_{n, 1}$ from both sides, our condition becomes same as \eqref{eq:re_explanation_nonparam} as used in \cite{tan2019doubly}. 
We now state the statistical error theorem: 
\begin{theorem}\label{thm:fixed_error_sparse}
Consider the function $\hat{\phi}$ defined by \eqref{eq:fixed_sparse_penalize}. Then, under the same assumptions as that of Theorem \ref{thm:fixed_design_wo_rsc} along with Assumption \ref{assn:rsc}, we have: 
\begin{equation}
    \|\hat f - f_0\|_n^2 = O_p\left(s_1 (\rho_{n, 1}^2+\lambda^2_{n,1})+ s_2 (\rho_{n, 2}^2+\lambda^2_{n,2})\right) \,.
\end{equation}
\end{theorem}
\noindent The proof of this theorem is presented in  Appendix \ref{sec:proof_fixed_high_dim}. The key difference between the rate obtained in Theorem \ref{thm:fixed_design_wo_rsc} and Theorem \ref{thm:fixed_error_sparse} is the absence of $\lambda_{n, i}$, which yields a faster rate under Assumption \ref{assn:rsc}. Similar discussion as that of after Theorem \ref{thm:fixed_design_wo_rsc} is in order; we need to choose $N_1, N_2$ to balance $\rho_{n, i}^2$ and $\lambda_{n, i}^2$. This leads to the following corollary: 
\begin{corollary}
    \label{cor:high_dim_rate_rsc}
    Choosing $N_1 = n^{1/(2(1 + 2\beta_1))}$ and $N_2 = n^{1/(2(1 + \beta_2))}$ (take the nearest integer if they are not integers), we obtain: 
$$
 \|\hat f - f_0\|_n^2 = O_p\left(s_1\left(n^{-\frac{2\beta_1}{1 + 2\beta_1}}\log^4{n} + \frac{\log{d}}{n}\right)  + s_2 \left(n^{-\frac{\beta_2}{1 + \beta_2}}\log^4{n} + \frac{\log{d}}{n}\right)\right) \,.
$$
\end{corollary}
In Section \ref{sec:minimax_high_dim}, we prove a matching lower bound to obtain that this rate is in fact minimax optimal up to log factors.

\subsection{Random design setting}
\label{sec:random_design_highdim}
\noindent
In the previous section, we have analyzed the fixed design model, i.e. $X_i$'s are assumed to be fixed. Here, we extend our analysis to the random design model, namely when $X_i$'s are random variables with distribution satisfying Assumption \ref{assn:density_bound}. 
As $X_i$'s are assumed to be random, we analyze the behavior of an estimator in terms of expected squared error loss, also known as generalization/out-of-sample error. Let us first highlight the difference between this subsection and the previous one. In the fixed design setup, the entire analysis hinges on the given $(X_1, \dots, X_n)$, i.e. we don't require to evaluate the performance of the predictor on some unobserved $X$. However, when we assume $X$'s are random, then it is imperative to evaluate the performance on the unseen $X$'s as there is a high probability (in fact probability is 1 if some component of $X$ is continuous) that in future we need to predict on new observations. Hence, in this setup we need to bound the expected squared error loss, i.e. $\|\hat f - f_0\|_{L_2(P_X)}^2$ instead of empirical squared error loss, i.e. $\|\hat f - f_0\|_n^2$. \par
A key step to bound this generalization error is to control the fluctuation of the corresponding empirical process, which ensures that a predictor, that performs well on the training samples, also performs well on the previously unobserved test samples. Typically, it is difficult to achieve such a guarantee using only the $L_1$-norm of the $L_2(\bbP_n)$ penalties as this penalty only controls the complexity of the functions on the observed samples. Although it might be possible to obtain a minimax optimal estimator by using $L_\infty$ penalty along with the $L_1$-norm of the $L_2(\bbP_n)$ penalties on the component functions, the optimization procedure becomes computationally challenging. To circumvent this issue, we propose a computationally efficient two-step procedure, which, at the cost of mildly stronger assumptions, allows us to estimate the true mean function optimally. We illustrate our estimator via Algorithm \ref{algo:random_design_algo}. 
\begin{algorithm}[!ht]
\DontPrintSemicolon
  \KwInput{Constant $c_1, c_2$ and penalty parameters $\lambda_{n, 1}$, $\lambda_{n, 2}$.}
  \KwOutput{$\hat f$: the estimator based on neural network.} 
   \KwData{Dataset $\{(X_1, Y_1), \dots, (X_n, Y_n)\}$.}
  \vspace{0.1in}
  Divide the dataset into two equal halves with $n/2$ data in each set (if $n$ is odd, then take $(n+1)/2$ data in the first half and $(n-1)/2$ in the second half). Denote the halves by $\cD_1$ and $\cD_2$. \\
  \vspace{0.1in}
  Using $\cD_1$, estimate the component functions $\{\hat \phi^{\text{init}}_j\}$ and $\{\hat \phi^{\text{init}}_{jk}\}$ by solving \eqref{eq:fixed_sparse_penalize}. \\
  \vspace{0.1in}
  Set $\hat S_1 = \{j: \|\hat \phi^{\text{init}}_j\|_n \ge c_1\lambda_{n, 1}\}$ and $\hat S_2 = \{(j < k): \|\hat \phi^{\text{init}}_{jk}\|_n \ge c_2 \lambda_{n, 2}\}$. \\
  \vspace{0.1in}
  Re-estimate component functions on the estimated active set $\hat S_1$ and $\hat S_2$ by minimizing (un-penalized) squared error loss:
  \begin{align}
  \label{eq:two_step_OLS_eq}
  & \{\hat \phi^{\text{final}}_j\}_{j \in \hat S_1}, \{\hat \phi^{\text{final}}_{jk}\}_{(j, k) \in \hat S_2} \notag \\
  & = \argmin_{\substack{\phi_j \in \cF_{NN, 1} \\ \phi_{j, k} \in \cF_{NN,2}}} \frac{1}{n}\sum_i \left(Y_i - \sum_{j \in \hat S_1}\phi_j(X_{ij}) - \sum _{(j, k) \in \hat S_2}\phi_{jk}(X_{ij}, X_{ik})\right)^2 \,.
  \end{align}  \\
  \vspace{0.1in}
  Return the final estimate $\hat f^{\text{final}} = \sum_{j \in \hat S_1}\hat \phi^{\text{final}}_j + \sum_{(j < k) \in \hat S_2}\hat \phi^{\text{final}}_{jk}$ \,.
\caption{Estimation under random design setting}
\label{algo:random_design_algo}
\end{algorithm} 
\par 
The key idea is as follows: we first split the samples into two (almost) equal halves. Based on the first half of the sample, we estimate the component functions by the same penalized procedure used in the fixed design setting (see\eqref{eq:fixed_sparse_penalize}). Next, we use hard thresholding on the estimators of univariate and bivariate component functions to estimate the active set (Step 3 of Algorithm \ref{algo:random_design_algo}). Denote the estimated active sets by $\hat S_1$ and $\hat S_2$ respectively. Note that our threshold levels are proportional to the penalty applied in \eqref{eq:fixed_sparse_penalize}. The constants $c_1, c_2$ will be mentioned explicitly in the proof. Once we have $\hat S_1, \hat S_2$, then we solve the unpenalized least square problem \textit{only on the active set} based on the second half of the data (Step 4 of Algorithm \ref{algo:random_design_algo}) to obtain the final estimate. \par 
Our algorithm is primarily motivated by SURE independent screening  \citep{fan2008sure, fan2010sure} and the idea of least square regression upon selecting an active set via LASSO or other penalized regression  \citep{belloni2013least}. SURE independent screening selects the active variables by performing marginal regression for 
additive models, whereas for sparse linear regression, \cite{belloni2013least} proposes to choose the active set using LASSO and then perform OLS on the selected subset to reduce bias. We combine these two ideas in Algorithm \ref{algo:random_design_algo}. First, we estimate the active subset of univariate and bivariate components $\hat S_1, \hat S_2$ (step 2-3 of Algorithm \ref{algo:random_design_algo}) and ensure that with high probability: i) $S_i \subseteq \hat S_i$ and ii) $|\hat S_i \cap S_i^c| = O(|S_i|)$. \par 
We need i) to avoid any potential bias that may occur by ignoring active components and ii) to guard against false positives, i.e. we do not select too many inactive variables. To ensure this, we need the following assumption:
\textcolor{black}{
\begin{assumption}
    \label{asm:signal_strength}
     Assume that all the component functions has minimal signal strength $r_n$, i.e. $\min_{j \in \cS_1} \|f_j^0\|_2 \ge r_n$, $\min_{(j,k) \in \cS_2} \|f_{j,k}^0\|_2 \ge r_n$ where: 
     $$
     r_n \gtrsim \sqrt{\left(s_1\left(n^{-\frac{2\beta_1}{1 + 2\beta_1}}\log^4{n} + \frac{\log{d}}{n}\right)  + s_2 \left(n^{-\frac{\beta_2}{1 + \beta_2}}\log^4{n} + \frac{\log{d}}{n}\right)\right)} \,.
     $$ 
\end{assumption}
}
\noindent Note that $r_n$ is precisely the rate of convergence obtained in Corollary \ref{cor:high_dim_rate_rsc}. Assumption \ref{asm:signal_strength} ensures that the signal in the active components is strong enough to be detected with high probability (ensuring the true active set is a subset of the estimated active set). Furthermore, by choosing an optimal penalty level $\lambda_{n, i}$ we control the number of false positives. 
Finally, in step 4 of Algorithm \ref{algo:random_design_algo}, we optimized squared error loss on the selected subset to produce the final estimate. The following theorem establishes rate of convergence of the estimator $\hat f_{\text{final}}$ obtained via Algorithm \ref{algo:random_design_algo}:  
\begin{theorem}
\label{thm:random_design_hard_threshold_ubd}
     Assume that the restricted strong convexity assumption (Assumption \ref{assn:rsc}) holds with high probability on $\{X_1, \dots, X_n\}$.  
    Then, under Assumption \ref{asm:signal_strength}, the estimator $\hat f^{\text{final}}$ satisfies, upto log-factors,
    $$
    \|\hat f^{\text{final}} - f_0\|_2^2 = O_p\left(s_1\lambda_{n, 1}^2 +  s_2\lambda_{n, 2}^2\right)\,.
    $$
    where $\lambda_{n, 1}$ and $\lambda_{n, 2}$ are same as in Theorem \ref{thm:fixed_error_sparse}. 
\end{theorem}
\noindent The proof of the Theorem will be deferred to Appendix \ref{sec:proof_ht}.

\subsection{Minimax lower bound}
\label{sec:minimax_high_dim}
In this section, we establish the minimax lower bound on the rate of convergence of the high dimensional two-way interaction model. As mentioned previously, our proofs can be extended to a general $k$-way interaction model. Our main theorem is as follows: 
\begin{theorem}
    \label{thm:minimax_high_dim}
    Suppose that we have $n$ observations from model \eqref{eq:define_sparse_model}. Then under Assumptions \ref{assn:identifiability}, \ref{assn:density_bound} we have: 
    $$
    \inf_{\hat f} \sup_{\substack{f \in \mathcal{F}_{\text{sp}} \\ X \sim P_X}} \|\hat f - f\|_2^2 \gtrsim \left(s_1\left(n^{-\frac{2\beta_1}{2\beta_1 + 1}} \vee \frac{\log{\left(d/s_1\right)}}{n} \right) + s_2\left(n^{-\frac{2\beta_2}{2\beta_2 + 2}} \vee \frac{\log{\left(d^2/s_2 \right)}}{n}  \right) \right)
    $$
\end{theorem}
\noindent It is immediate from Theorem \ref{thm:minimax_high_dim} and Theorem \ref{thm:random_design_hard_threshold_ubd} that our estimator of the mean function based on neural networks is indeed minimax optimal up to log factors.  \par 
Although the proof techniques of Theorem \ref{thm:minimax_high_dim} and Theorem \ref{thm:minimax_lower_bound_growing_d} are similar, there are some important changes in the construction for the alternatives since  we need to enforce sparsity here. We would also like to point out that our proof is different from that of \cite{raskutti2012minimax} since we do not assume that the components of the mean function belong to RKHS. Here we sketch the main idea of the proof briefly: each univariate component function is assumed to be $\beta_1$-smooth and each bivariate component function is assumed to be $\beta_2$-smooth. We first construct a collection of alternatives for these components and then we take a sparse combination of them. Finally, constructing alternatives along with a proper choice of the relevant hyper-parameters yield the lower bound. The proof is deferred to the Appendix.

\section{Conclusion}\label{sec:conclusion}

Deep neural networks have achieved tremendous success in nonparametric function estimation. Yet due to the intrinsic difficulty of ``curse-of-dimensionality'', they can only handle nonparametric functions of low-dimension, without excessive restrictions on the function classes. This calls for low-dimensional nonparametric interaction models. At the same time, modern big data applications often involve nonparametric regression with a large number of predictors.  Yet, most statistical theories on neural networks focus only on finite-dimensional regression.  These give rise to the imminent need for the study of nonparametric interaction models in diverging dimensions and understanding the impact of dimensionality in such structured nonparametric models.

This paper contributes critically to understanding the performance of neural networks in low-order interaction models in diverging dimensions.  An important conclusion of our study is that estimated components should have low biases in order to avoid unnecessary error accumulations and this is achieved by our newly introduced debiasing techniques.  
For slowly diverging dimensional problems, no additional regularization is needed.  With proper debiasing,direct least-squares estimation on the structured neural networks is shown to achieve a rate of convergence that matches with a newly established minimax low bound, namely, it is optimal.  In a high-dimensional setting, sparsity assumption on the interaction terms is necessary.  We appeal to the penalized least-squares estimation and screening techniques (for random design) and show that the resulting procedure is minimax optimal by establishing a matching lower bound.  Our results provide a comprehensive view of the performance of neural networks for the structured nonparametric model in diverging dimensions, which are critical to modern data science and necessary for neural networks to succeed.

\section{Proof of Theorem \ref{thm:approx_growing_d}}
\label{sec:proof_approx_growing_d}
Suppose $f:[0, 1]^d \to \reals$ is $\beta$-smooth, i.e. it is $\lfloor \beta \rfloor$ times differentiable with bounded derivatives. Then by Taylor series approximation around some point $x_0$, we have: 
\begin{align*}
f(\bx) & = \sum_{|\balpha| \le \beta}\frac{1}{\balpha !}\frac{\partial^\balpha f(\bx_0)}{\partial x^\balpha}(\bx - \bx_0)^\balpha + \sum_{|\balpha| = \beta} \frac{1}{\balpha !} \frac{\partial^\balpha f(\bx_0 + \lambda_\balpha \bx)}{\partial x^\balpha}(\bx - \bx_0)^\balpha \\
& := T(\bx) + R(\bx) 
\end{align*}
where we use standard multi-dimensional Taylor series notation: for any $\balpha = (\balpha_1, \dots, \balpha_d)$, set $|\balpha| = \sum_{j=1}^{d} \balpha_j$. Let $\partial^\balpha/\partial x^\balpha$ denotes $\partial^{\balpha_1 + \dots + \balpha_d}/\partial x_1^{\balpha_1} \partial x_2^{\balpha_2} \dots \partial x_d^{\balpha_d}$. For any vector $\bv$, let $\bv^\balpha$ denote $v_1^{\balpha_1}\dots v_d^{\balpha_d}$. The basic idea of approximating a smooth function by a neural network involves the following three key steps: 
\begin{enumerate}
    \item Divide the domain $[0, 1]^d$ into small grids. 
    \item On each grid, perform Taylor series expansion around some fixed point, say the midpoint of the corner point of the grid. 
    \item Approximate the Taylor polynomials by neural networks. 
\end{enumerate}
\begin{figure}
    \centering
    \includegraphics[scale=0.3]{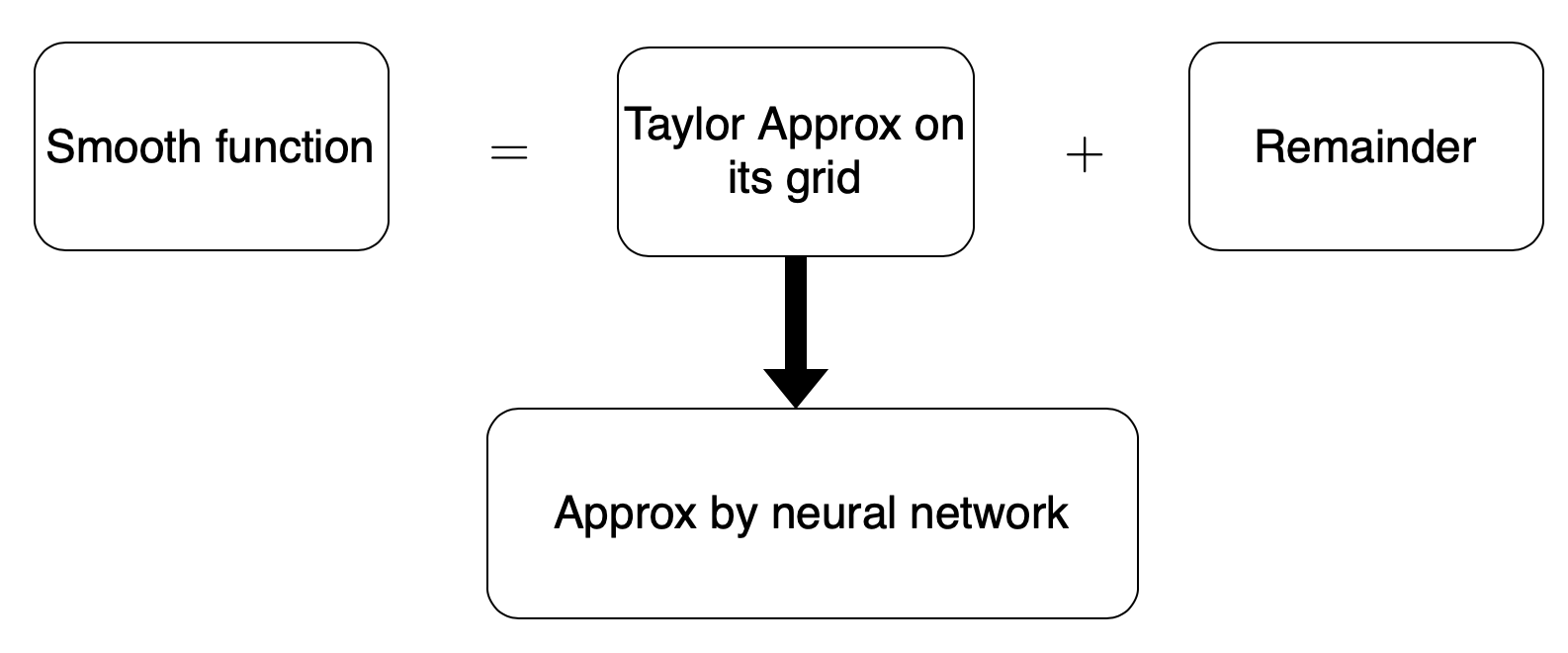}
    \caption{Approximating smooth function by neural network}
    \label{fig:taylor_approx}
\end{figure}
Figure \ref{fig:taylor_approx} presents a flowchart with these three steps. Following the notations of \cite{lu2021deep}, the best neural network approximator can be written as: 
\begin{equation}
\label{eq:nn_approx_main}
\phi(x) = \sum_{|\alpha| \le \beta} \varphi\left(\frac{\phi_\alpha\left(\bPsi(x)\right)}{\alpha!}, \bP_\alpha(x - \bPsi(x))\right)
\end{equation}
where, $\bPsi(x)$ maps the point $x$ to the corner point of its grid, $\phi_\alpha$ approximates the value of the derivatives, $\bP_\alpha$ approximates the polynomial $\bx^\balpha$ and finally $\varphi$ (the outer function) approximates the operation $xy$ by $\phi(x, y)$. Next we describe how to perform addition using neural networks. Given $k$ neural networks each having width $N$ and depth $L$, we can add them the following way: first send $x$ to $2k$ neural networks $\{\phi_1, \phi_2, \dots, \phi_k\}$ and $\{-\phi_1, \dots, -\phi_k\}$. Then the final output is $\sum_{j=1}^k \left(\sigma(\phi_j(x)) - \sigma(-\phi_j(x))\right)$. As $\sigma(\cdot)$ is ReLU, we have the identity $\sigma(x) - \sigma(-x) = x$. See Figure \ref{fig:deepnet_addition} for a visual description. 
Therefore, the sum can be performed using a NN with width $2kN$ and depth $L+1$. If $k$ is large, then one can perform the same operation via a neural network of width $kN \vee 2k$ and width $L+2$. In \eqref{eq:nn_approx_main} the value of $k$ is of the order $\beta^d$ which is fixed, so the order of width and depth remains unchanged. \par
The neural network approximation, as presented in \eqref{eq:nn_approx_main}, has five sources of error (Table \ref{tab:error_approx}): i) $E_1$, the error of approximating $f$ by its Taylor approximation on a grid, ii) $E_2$, the error of approximation of $\bPsi(\bx)$, iii) $E_3$, the error of approximation of $\phi_\alpha$, iv) $E_4$, the error of approximating $\bx^\balpha$ by $\bP_\alpha$ and finally v) $E_5$, the error of approximating $xy$ by $\phi(x, y)$. For any given $N$ and $L$, suppose we divide $[0, 1]$ into $K$ equispaced intervals with length $1/K$ where $K = \lfloor N^{1/d} \rfloor^2 \lfloor L^{2/d}  \rfloor$. If we use a neural network of width $C_1 N \log{N}$ and width $C_2 L \log{L}$ for some $C_1,C_2>0$, then we commit the following approximation error: 
\begin{table}
\centering
\label{tab:error_approx}
\begin{tabular}{||c c c c||} 
 \hline
 Errors & & & Order \\ [0.7ex] 
 \hline\hline
 $E_1$ & & & $K^{-\beta} \sim (NL)^{-2\beta/d}$\\ 
 \hline
 $E_2$ & & & $0$ \\
 \hline
 $E_3$ & & & $2(NL)^{-2\beta}$\\
 \hline
 $E_4$ & & & $9\beta (N+1)^{-7\beta L}$\\
 \hline
 $E_5$ & & & $216(N+1)^{-7\beta L}$\\
 \hline
\end{tabular}
\caption{Approximation error of neural network (cf. \cite[Table 2]{lu2021deep})}
\end{table}
Therefore for $\beta \ge 1$ and $d \ge 1$, the only term where the effect of dimension creeps in is the Taylor approximation error which is $E_1$. The NN $\bPsi$ is a step function, which $\bx$ to the corner point of its grid without any error. $\phi_\alpha$ is a point fitting function, which maps the corner point of the grid to its $\balpha^{th}$ derivative of $f$ and commits error of the order $(NL)^{-2\beta}$ which does not depend on its dimension. $\bP_\balpha$ is a neural network which approximates $\bx^{\balpha}$ and has error of the order $(N+1)^{-7\beta L}$ and so is the neural network $\varphi$ which approximates the bivariate function $(x, y) \rightarrow xy$. 
\\\\
\noindent
\hbox{\bf Step 1: (Approximating one dimensional component)} For each one dimensional component $f_j$, there exists a neural network $\phi_j$ (by the above mentioned construction) with width $C_1 N_1\log{N_1}$ and depth $C_2 L_1\log{L_1}$ such that, for some $C>0$,
$$
\|f_j - \phi_j\|_\infty \le C(N_1L_1)^{-2\beta_1} \,.
$$
Define $I_j = \int_0^1 \phi_j(x) \ dx$ and another neural network $\widetilde \phi_j(x) = \phi_j(x) - I_j$. As subtracting a constant tantamounts to changing the bias of the last layer, it does not change the architecture. Therefore, the width and depth of $\widetilde \phi_j$ is same as $\phi_j$. Furthermore we have: 
$$
\|f_j - \widetilde \phi_j\|_\infty = \left\|f_j - \int_0^1 f_j - \phi_j - \int_0^1 \phi_j\right\|_\infty \le 2\|f_j - \phi_j\|_\infty \le 2C(N_1L_1)^{-2\beta_1} \,.
$$
{\bf Step 2: (Approximating two dimensional component) }Now consider a two dimensional function, say $f_{12}$. As before, following \cite{lu2021deep}, we can construct a neural network $\phi_{12}$ with width $O(N_2\log{N_2})$ and depth $O(L_2\log{L_2})$ such that
\begin{equation}
\label{eq:two_dim_approx_1}
\|f_{12} - \phi_{12}\|_\infty \le C(N_2L_2)^{-\beta_2}. 
\end{equation}
However, by \eqref{eq:nn_approx_main}, $\phi_{12}$ does not approx $f_{12}$ directly, it only approximates the Taylor expansion $T_{12}(\bx)$. Hence, there exists some $L^\star(\beta)>0$ such that for any $L \ge L^\star(\beta)$, 
\begin{equation}
\label{eq:two_dim_approx_2}
\|f_{12}(\bx) - \phi_{12}(\bx)\|_{\infty} \le \underbrace{\|T_{12}(\bx) - \phi_{12}(\bx)\|_\infty}_{O((N_1L_1)^{-2\beta_2})} + \underbrace{\|R_{12}(\bx)\|_\infty}_{O((N_2L_2)^{-\beta_2})}.
\end{equation}
Note that $T_{12}(\bx)$ is a bivariate polynomial. Although the overall approximation error of the neural network is $O((N_2L_2)^{-\beta_2})$, the approximation error of $T_{12}$ by $\phi_{12}$ is much faster and can be assumed to be faster than $(N_2L_2)^{-2\beta_2}$ if $L \ge L^\star(\beta_2)$ for some constant $L^\star(\beta_2)$ only depending on $\beta_2$. If we define $T_1(x)$ (resp. $T_2$(y)) as $T_1(x) = \int_0^1 T_{12}(x, y) \ dy$ (resp. $T_2(y) = \int_0^1 T_{12}(x, y) \ dx$), then both $T_1$ and $T_2$ are univariate polynomial of degree $\le \beta_2-1$ and hence $\beta_2$ smooth. Consequently, we can construct NN $\xi_1$ and $\xi_2$ with width of $O(N_2\log{N_2})$ and depth of $O(L_2\log{L_2})$ such that: 
\begin{align}
\label{eq:one_dim_marginal_approx}
\|T_1(x) - \xi_1(x)\|_\infty & \le C(N_2L_2)^{-2\beta_2} \,, \\
\|T_2(y) - \xi_2(y)\|_\infty & \le C(N_2L_2)^{-2\beta_2} \,.
\end{align}
We define  the constant $I_{12} = \iint \phi_{12}(x, y) \ dx \ dy$. Finally, we define the estimator by $\widetilde \phi_{12}(x, y) := \phi(x, y) - \xi_1(x) - \xi_2(y) + I_{12}$. Let us comment about the architecture of $\tilde \phi_{12}$. As mentioned before, adding $I_{12}$ does not change the architecture. Also as $\phi_{12}, \xi_1, \xi_2$ has width $O(N_2\log{N_2})$ and depth $O(L_2\log{L_2})$, $\widetilde \phi_{12}$ also has width $O(N_2\log{N_2})$ and depth $O(L_2\log{L_2})$ (see Figure \ref{fig:deepnet_addition}). 
Next, we show that the $L_\infty$ distance of $\widetilde \phi_{12}$ form $f_{12}$ remains of the order $(N_2L_2)^{-\beta_2}$. Recall that the marginals of $f_{12}$ are $0$. Using triangle inequality we have:  
\allowdisplaybreaks
\begin{align*}
    \|f_{12} - \widetilde \phi_{12}\|_{\infty} & = \|f_{12} - \phi_{12} + \xi_1 + \xi_2 - I_{12}\|_{\infty} \\
    & \le \|f_{12} - \phi_{12}\|_{\infty} + \left\|\int_0^1 f_{12} \ dy - \xi_1\right\|_\infty \\
    & \qquad \qquad \qquad + \left\|\int_0^1 f_{12} \ dx - \xi_2\right\|_\infty + \left\| \iint_{[0,1]^2} f_{12} \ dx \ dy - I_{12}\right\|_\infty \\
    & \le 2\|f_{12} - \phi_{12}\|_{\infty} + \left\|T_1 + \int_0^1 R_{12} \ dy - \xi_1\right\|_\infty \\
    & \qquad \qquad \qquad + \left\|T_2 + \int_0^1 R_{12} \ dx - \xi_2\right\|_\infty \hspace{0.1in} \Big[\text{As }I_{12} = \iint \phi_{12}\Big] \\
    & \le 2\|f_{12} - \phi_{12}\|_{\infty} + \|T_1 - \xi_1\|_\infty + \|T_2 - \xi_2\|_\infty \\
    & \qquad \qquad + \left\| \int_0^1 R_{12} \ dy\right\|_\infty + \left\|\int_0^1 R_{12} \ dx\right\|_\infty \\
    \vspace{0.1in}
    & \le 2\|f_{12} - \phi_{12}\|_{\infty} + \|T_1 - \xi_1\|_\infty + \|T_2 - \xi_2\|_\infty + 2\|R_{12}\|_\infty \\
    \vspace{0.1in}
    & \le C(N_2L_2)^{-\beta_2} + 2C(N_2L_2)^{-2\beta_2} + 2C(N_2L_2)^{-\beta_2} = 5C(N_2L_2)^{-\beta_2} \,.
\end{align*}
Here the first error follows from equation \eqref{eq:two_dim_approx_1}, second and third from equation \eqref{eq:one_dim_marginal_approx} and the last one from equation \eqref{eq:two_dim_approx_2}. 
\\\\
\noindent
We next show that $\left\|\int_0^1 \left(\widetilde \phi_{12}- f_{12} \right) \ dy\right\|_\infty \lesssim (NL)^{-2\beta}$. Since $\int_0^1 f_{12} \ dy =0$, the claim is true since
\allowdisplaybreaks
\begin{align*}
    \left\|\int_0^1 \widetilde \phi_{12} \ dy\right\|_\infty & = \left\|\int_0^1 \phi_{12} \ dy - \xi_1 - \int_0^1 \xi_2 \ dy + I_{12} \right\|_\infty \\
    & \le \left\|\int_0^1 (\phi_{12} - T_{12})\ dy + T_1 - \xi_1\right\|_\infty \\
    & \qquad \qquad + \left\|\int_0^1 (\xi_2 - T_2) \ dy + \int_0^1 T_2 \ dy - I_{12}\right\|_\infty \\
    & \le \|\phi_{12} - T_{12}\|_{\infty} + \|\xi_1 - T_1\|_\infty + \|\xi_2 - T_2\|_\infty\\
    & \qquad + \left\|\iint_{[0,1]^2} T_{12} \ dx \ dy -  \iint_{[0,1]^2} \phi \ dx \ dy\right\|_\infty \\
    & \le  2\|\phi_{12} - T_{12}\|_{\infty} + \|\xi_1 - T_1\|_\infty + \|\xi_2 - T_2\|_\infty
    \le 4C(N_2L_2)^{-2\beta_2} 
\end{align*}
Here, the first inequality follows from equation \eqref{eq:two_dim_approx_2} second and third inequality follows form \eqref{eq:one_dim_marginal_approx}. 
\\\\
\noindent
{\bf Step 3: (Combining all the bounds)} Here, we combine the bounds for one dimensional and two dimensional components to prove Theorem \ref{thm:approx_growing_d}. Define the function $\widetilde \phi$ as: 
$$
\widetilde \phi = \sum_{j=1}^d \widetilde \phi_j + \sum_{k < l}\widetilde \phi_{kl} \,.
$$
Ging back to our argument for summing multiple neural network, the width of $\sum_j \tilde \phi_j$
is $O(dN_1\log{N_1})$ and depth is $O(L_1\log{L_1})$, the width of $\sum_{k < l}\widetilde \phi_{kl}$ 
 is $O((d(d-1)/2)N_2\log{N_2})$ and depth is $O(L_2\log{L_2})$. Furthermore, we have:
\allowdisplaybreaks
\begin{align*}
    & \bbE\left[\left(f_0(X) - \widetilde \phi(X)\right)^2\right] \le p_{\max} \int_{[0, 1]^d}\left(f_0(\bx) - \widetilde \phi(\bx)\right)^2 \ d\bx \\
    & = p_{\max}\int_{[0, 1]^d}\left(\sum_j f_{0, j}(x_j) + \sum_{i < j} f_{0, ij}(x_i, x_j) - \sum_j \widetilde  \phi_j(x_j) - \sum_{i < j}\widetilde \phi_{ij}(x_i, x_j)\right)^2 \ d\bx  \\
    & \le 2p_{\max}\int_{[0, 1]^d}\left(\sum_j f_{0, j}(x_j) - \sum_j \widetilde  \phi_j(x_j)\right)^2 \ d\bx + 2p_{\max}\int_{[0, 1]^d}\left(\sum_{i < j} f_{0, ij}(x_i, x_j) - \sum_{i < j}\widetilde \phi_{ij}(x_i, x_j)\right)^2 \ d\bx \\
    & = 2p_{\max}\sum_{j=1}^d \int_0^1 \left(f_{0, j}(x_j) - \widetilde \phi_j (x_j)\right)^2 \ dx_j + 2p_{\max}\sum_{i < j}\int_{[0, 1]^2}\left(f_{0, ij}(x_i, x_j) - \widetilde \phi_{ij} (x_i, x_j)\right)^2 \ dx_i \ dx_j \\
    &+ 2p_{\max}\sum_{j \neq j'}\int_0^1 \left(f_{0, j}(x_j) - \widetilde \phi_j (x_j)\right) \ dx_j\int_0^1 \left(f_{0, j'}(x_{j'}) - \widetilde \phi_{j'} (x_{j'})\right)^2 \ dx_{j'} \\
    &+2p_{\max}\int_{[0, 1]^4}\sum_{(i, j) \neq (k, l)} \left(f_{0, ij}(x_i, x_j) - \widetilde \phi_{ij} (x_i, x_j)\right)\left(f_{0, kl}(x_k, x_l) - \widetilde \phi_{kl} (x_k, x_l)\right) \ dx_idx_jdx_kdx_l \\
    & \le 2p_{\max}\left[\sum_{j=1}^{d} \|f_{0, j} - \widetilde \phi_j\|^2_\infty + \sum_{i < j} \|f_{0, ij} - \widetilde \phi_{ij}\|_\infty^2\right] \\
    & \qquad +2p_{\max}\sum_{\substack{i \neq (k, l) \\ j \neq (k, l)}} \int_{[0, 1]^2} \widetilde \phi_{ij} \ dx_i \ dx_j \int_{[0, 1]^2} \widetilde \phi_{kl} \ dx_k \ dx_l + \\
    & \qquad \qquad + 2p_{\max}\sum_{i, j , k} \int_0^1 \left(\int_0^1 \widetilde \phi_{ij} \ dx_j\right)\left(\int_0^1 \widetilde \phi_{ik} \ dx_k\right) \ dx_i\\
    & \le 2C^2p_{\max}\left[d(N_1L_1)^{-4\beta_1} + \dbinom{d}{2} (N_2L_2)^{-2\beta} + \dbinom{d}{4}(N_2L_2)^{-4\beta} + \dbinom{d}{3} (N_2L_2)^{-4\beta}\right] \\
    & \le C_3 \left(d(N_1L_1)^{-4\beta_1} + \binom{d}{2}(N_2L_2)^{-2\beta_2}\right), \,
\end{align*}
for some $C, C_3>0$ since $d^2(N_2L_2)^{-2\beta_2} \rightarrow 0$ by Assumption \ref{asm:d_growth}. Here the penultimate inequality follows by combining the results of Step 1 and Step 2. This completes the proof.

\begin{funding}
This paper is supported by ONR  N00014-22-1-2340 and the NSF grants DMS-2052926, DMS-2053832, DMS-2210833.
\end{funding}

\bibliography{Refs}
\bibliographystyle{ims}
\appendix
\newpage

\section{Proofs}\label{sec:proofs}

\subsection{Proof of Theorem \ref{thm:statistical_error_lowdim}}
\label{sec:proof_statistical_error_lowdim}
We want to invoke \cite[Theorem 3.2.5]{vaart1997weak}. For any $\phi = \phi_1 + \phi_2$ (where $\phi_1 \in \cF_{NN}^1, \phi_2 \in \cF^2_{NN}$), 
define the population risk function $R(\phi)$ (resp. empirical risk $\hat R_n(\phi)$) as 
    $$ 
    R(\phi) = \bbE[(Y - \phi(X))^2] \ \ \ \left(\text{resp. } \hat R(\phi) = \frac1n \sum_i (Y_i - \phi(X_i))^2\right) 
    $$ 
    Since that $\eps$ is independent of $X$ with mean $0$, it is immediate that: 
    $$
    R(\phi) - R(f_0) = \bbE[(f(X) - f_0(X))^2] \triangleq d^2(\phi, f_0) \,.
    $$
Define $\cF_n := \{\phi = \phi_1 + \phi_2: \phi_1 \in \cF^1_{NN}, \phi_2 \in \cF^2_{NN}\}$ (i.e. $\cF_n = \cF_{NN}^1 + \cF_{NN}^2$) and $\omega_n = d(\phi^\star, f_0)$, i.e. the approximation error of the class of neural network. It is immediate that: 
$$
\textstyle
\cF_n \subseteq \cF_{NN}\left(d, c_1 \left(dN_1\log{N_1} + \binom{d}{2}N_2\log{N_2}\right), c_2, W_n, 1\right)
$$
with 
$$
W_n = c_3\left(d(N_1\log{N_1})^2 + \binom{d}{2}\left(N_2\log{N_2}\right)^2\right) \,,
$$
for some constant $c_1, c_2, c_3$ independent of $(n, d)$. Therefore, VC dimension of $\cF_n$ is less than or equal to $V_n \triangleq c_4 W_n \log{W_n}$ for some $c_4>0$. We use this in the subsequent analysis. Let $\zeta_n(\delta)$ be such that: 
    \begin{equation}
    \label{eq:max_ineq_1}
    \bbE\left[\sup_{\phi \in \cF_n: d(\phi, f_0) \le \delta} \left|\left(R_n(\phi) - R(\phi)\right) - \left(R_n(f_0) - R(f_0)\right)\right|\right] \lesssim \frac{\zeta_n(\delta)}{\sqrt{n}} \,.
    \end{equation}
This function $\zeta_n$ is called the modulus of continuity, which is used to bound the fluctuation of the empirical process around the population process in a neighborhood of the true function. $\zeta_n(\delta)$ intrinsically depends on the complexity of the underlying function class. We quantify this dependency through the following maximal inequality (e.g. Theorem 5.2 of \cite{chernozhukov2014gaussian}), which we state here for the convenience of the readers:

\begin{lemma}[Theorem 5.2 of \cite{chernozhukov2014gaussian}]
\label{lem:max_ineq}
    Define the operator $\bbG_n = \sqrt{n}(\bbP_n - P)$. Consider a collection of functions $\cF$ with envelope function $F$. Assume that $F \in L_2(P)$. Then: 
$$
\bbE\left[\sup_{f \in \cF} \left| \bbG_n(f)\right|\right] \lesssim J(r, \cF, F)\|F\|_2 + \frac{\|M\|_2 J^2(r, \cF, F)}{r^2 \sqrt{n}}
$$
for $r = \sigma/\|F\|_2$,  where
\begin{align*}
    J(\tau, F, \cF) & = \int_0^\tau \sup_Q \sqrt{1 + \log{N\left(\eps \|F\|_2, \cF, L_2(Q)\right)}} d \eps \\
    M & = \max_{1 \le i \le n} F(X_i), \qquad  \sigma^2 = \sup_{f \in \cF} Pf^2   \,.
\end{align*}
\end{lemma}

\noindent We use Lemma \ref{lem:max_ineq} to find the function $\zeta_n(\delta)$. 
For any $\phi$, define a function $g \equiv g(\phi)$ on the space of $(X, \eps)$ as: 
\begin{align}
    g(X, \eps) & = (f_0(X) - \phi(X) + \eps)^2 - \eps^2 \notag \\
\label{eq:g_val} & = (\phi(X) - f_0(X))^2 + 2(\phi(X) - f(X))\eps. 
\end{align}
Define the collection of $\cG$ as $\cG = \left\{g(\phi) : \phi \in \cF_{n} \right\}$. This definition implies:  
$$
\bbE\left[\sup_{\phi: d(\phi, f_0) \le \delta} \left|\left(R_n(\phi) - R(\phi)\right) - \left(R_n(f_0) - R(f_0)\right)\right|\right] = \frac{1}{\sqrt{n}}\bbE\left[\sup_{g \in \cG_\delta} \left|\bbG_n(g)\right|\right].
$$
Here we $\cG_\delta$ as the collection of all $g(\phi)$'s, such that $d(\phi, f_0) \le \delta$. 
To apply Lemma \ref{lem:max_ineq}, we need to i) quantify $\sigma^2$, ii) find an envelope of $\cG_\delta$ and iii) bound the logarithm of covering number. For i) observe that: 
\begin{align*}
    \sup_{g \in \cG: d(f, f_0) \le \delta} \bbE[g(X)^2] \le \delta^2(2 + 8\sigma_\eps^2)  \triangleq C_{\sigma_\eps}\delta^2 \,.
\end{align*}
Next, to construct an envelope, note that \eqref{eq:g_val} implies: 
$$
|g(X, \eps)| \le 9B^2 + 6B\eps \,.
$$
for all $(X, \eps)$. This follows from the fact that $\|\phi\|_\infty \le 2B$ and $\|f_0\|_\infty \le B$. Therefore, we can take the envelope $G$ as $G(X, \eps) = 9B^2 + 6B\eps$. 
To bound the logarithm of the covering number of $\cG_\delta$,  
we use Lemma 9.9 of \cite{kosorok2008introduction}. Note that, each $g=g_1 + g_2 + h g_3$ where $g_1 = (\phi - f_0)^2_+$, $g_2 = (\phi - f_0)^2_-$, $g_3 = 2(\phi - f_0)$ and $h(x, \eps) = \eps$. We calculate the VC dimensions of each class of functions: 
\begin{itemize}
    \item The VC dimension of $\cG_1 = \left\{g_1(\phi): \phi \in \cF_{NN}\right\}$ is $\le V_n$. To see this, write $g_1(x, \eps) = m_1 \circ (\phi - f_0) \circ m_2(x, \eps)$ where $m_2(x, \eps) = x$ and $m_1(x) = x_+^2$. As precomposing by a fixed function $\phi$, subtracting a fixed function $f_0$ and post-composing by a monotone function does not change the VC dimension, the claim follows. Similar argument establishes that VC dimension of $\cG_2 = \left\{g_2(f): f \in \cF\right\}$ is $\le V_n$. 

    \item VC dimension of $\cG_3$, the collection of $g_3h$ is $\lesssim V_n$. To see this, note that VC dimension of $g_3$ is $\le V_n$. As before, $g_3 = 2 (\phi - f_0) \circ \psi$ as we have already mentioned that pre-composing by $\psi$, subtracting $f_0$ and multiplying by by a fixed function (here $2h$) does not change the VC dimension.
\end{itemize}
Therefore, we have established that VC-dim of $\cG_i$ is $\lesssim V_n$. This implies, via Hausslers's bound \cite[Theorem 2.6.7]{van1996weak}): 
$$
\sup_Q \log{N\left(\eps \|G_i\|_2, \cG_i, L_2(Q)\right)} \lesssim
\log{KV_n} + 2V_n\log{\left(\frac{2\sqrt{2e}}{\eps}\right)} \,.
$$
As $\cG_\delta \subseteq \cG_1 + \cG_2 + \cG_3$, we have: 
\begin{align*}
   \log{N\left(\eps \|G\|_2, \cG_\delta, L_2(Q)\right)} & \lesssim \sum_{j = 1}^3 \log{N\left(\eps \|G_i\|_2, \cG_i, L_2(Q)\right)} \\
   & \le 3\log{KV_n} + 6V_n\log{\left(\frac{2\sqrt{2e}}{\eps}\right)} \\
   & \le 9V_n\log{\left(\frac{2\sqrt{2e}}{\eps}\right)} \,,
\end{align*}
where the last inequality holds as soon as $\log{KV_n} \le V_n$ and $\eps \le 2\sqrt{2/e}$. Therefore, we can use the last inequality for all small $\eps$ and for all large enough $V_n$. Taking $r = \sqrt{C_{\sigma_\eps}}\delta/\|G\|_2 $ we have: 
\begin{align*}
    J(r, \cG_\delta, G) & \le \int_0^r \sqrt{1 + 9V_n\log{\left(\frac{2\sqrt{2e}}{\eps}\right)}} \ d\eps \\
    & \le 4\sqrt{V_n}\int_0^r \sqrt{\log{\left(\frac{2\sqrt{2e}}{\eps}\right)}} \ d\eps \\
    & = 8\sqrt{2eV_n}\int_0^{r/2\sqrt{2e}} \sqrt{\log{\left(\frac{1}{\eps}\right)}} \ d\eps \\
    & \le 4r\sqrt{V_n}\sqrt{\log{\left(\frac{2\sqrt{2e}}{r}\right)}}.
\end{align*}
Here the second inequality holds as long as $V_n\log{(2\sqrt{2}/\eps)} \ge 1/7$, i.e. $\eps < 2\sqrt{2}\exp{\left(-1/(7V_n)\right)}$. This holds for all large $n$ as long as $r < 1$. Hence we conclude via Lemma \ref{lem:max_ineq}: 
\begin{align*}
    \bbE\left[\sup_{g \in \cG_\delta} \left|\bbG_n(g)\right|\right] & \lesssim  \delta \sqrt{V_n\log{\left(\frac{1}{\delta}\right)}} + \frac{V_n}{\sqrt{n}}\log{\left(\frac{1}{\delta}\right)}\sqrt{\bbE\left[\max_i (1 + 2|\eps_i|)^2\right]}  \\
    & \lesssim  \delta \sqrt{V_n\log{\left(\frac{1}{\delta}\right)}} + \frac{V_n \sqrt{\log{n}}}{\sqrt{n}}\log{\left(\frac{1}{\delta}\right)} \triangleq \phi_n(\delta) \,.
    \end{align*}
Now we establish the rate of convergence. Set $r_n = \left(\sqrt{\frac{V_n}{n}\log{\frac{n}{V_n} }\sqrt{\log{n}}} + \omega_n\right)^{-1}$, where $\omega_n = d(\phi^\star, f_0)$ is the approximation error defined before. Note that $r_n$ satisfies the equation $\phi_n(r_n^{-1}) \lesssim r_n^{-2}\sqrt{n}$ since
\begin{align*}
    r_n^2\phi_n\left(\frac{1}{r_n}\right) & \lesssim r_n\sqrt{V_n \log{r_n}} + r_n^2 \frac{V_n \sqrt{\log{n}}}{\sqrt{n}}\log{r_n} \\
    & \le \sqrt{n}\left[\left(\frac12\frac{\log{\frac{n}{V_n} - \log{\log{\frac{n}{V_n}}} - \frac12 \log{\log{n}}}}{\log{\frac{n}{V_n} }\sqrt{\log{n}}}\right)^{1/2} + \left(\frac12\frac{\log{\frac{n}{V_n} - \log{\log{\frac{n}{V_n}}} - \frac12 \log{\log{n}}}}{\log{\frac{n}{V_n} }}\right)\right] \\
    & \lesssim \sqrt{n}
\end{align*}

We now use shelling argument to establish the rate of convergence. Fix $t > 1$, and define $A_j := \{\cF_{NN} : 2^{j-1}t \le r_n d(f_0, \phi) \le 2^jt\}$.  It follows from the definition of ERM, we have
\allowdisplaybreaks
\begin{align}
& \bbP\left(r_n d(\hat \phi, f_0) \ge t\right) \notag = \bbP\left(\sup_{\substack{\phi \in \cF_{NN} \\ r_n d(f_0, \phi) \ge t}} \hat R_n(\phi^\star) - \hat R_n(\phi) \ge 0 \right) \notag \\
& \le \sum_{j=1}^\infty \bbP\left(\sup_{\phi \in A_j} \hat R_n(\phi^\star) - \hat R_n(\phi) \ge 0 \right) \notag \\
 & \le \sum_{j=1}^\infty \bbP\left(\sup_{\phi \in A_j} \left(R_n(\phi^\star) - R(\phi^\star)\right) - \left(R_n(\phi) - R(\phi)\right) \ge \inf_{\phi \in A_j} R(\phi) - R(\phi^\star) \right) \notag \\
\label{eq:rate_1} & \le \sum_{j=1}^\infty \frac{\bbE\left[\sup_{\phi \in A_j} \left|\left(R_n(\phi^\star) - R(\phi^\star)\right) - \left(R_n(\phi) - R(\phi)\right)\right|\right]}{\inf_{\phi \in A_j} R(\phi) - R(\phi^\star) }
\end{align}
We bound the numerator and denominator separately. For the denominator note that: 
\begin{align*}
    \inf_{f \in A_j} R(\phi) - R(\phi^\star)  & \ge \inf_{f: d_{f, f_0} \ge 2^{j-1}tr_n^{-1}} R(\phi) - R(f_0) + R(f_0) - R(\phi^\star) \\
    & \ge  2^{2j-2}t^2r_n^{-2} - \omega^2_n 
\end{align*}
For the numerator: 
\begin{align*}
    & \bbE\left[\sup_{\phi \in A_j} \left|\left(R_n(\phi^\star) - R(\phi^\star)\right) - \left(R_n(\phi) - R(\phi)\right)\right|\right] \\
    & \le \bbE\left[\sup_{f: d(f, f_0) \le 2^j tr_n^{-1}} \left|\left(R_n(\phi^\star) - R(\phi^\star)\right) - \left(R_n(\phi) - R(\phi)\right)\right|\right] \lesssim \frac{\zeta_n\left(2^j tr_n^{-1}\right)}{\sqrt{n}} \,.
\end{align*}
Here, the second last equation holds because $d(\phi^\star, f_0) = \omega_n \le 2^j tr_n^{-1}$ as $r_n^{-1} \ge \omega_n$. Putting this bound in equation \eqref{eq:rate_1} we have: 
\begin{align*}
    \bbP\left(r_n d(\hat f, f_0) \ge t\right) & \le \sum_{j=1}^\infty \frac{\phi_n\left(2^j tr_n^{-1}\right)}{\sqrt{n}\left(2^{2j-2}t^2r_n^{-2} - \omega_n^2\right) } \\
    & \le \sum_{j=1}^\infty \frac{2^j t\phi_n\left(r_n^{-1}\right)}{\sqrt{n}\left(2^{2j-2}t^2r_n^{-2} - \omega^2_n \right) } =  \sum_{j=1}^\infty \frac{2^j t r_n^2 \phi_n(r_n^{-1})n^{-1/2}}{2^{2j-2}t^2 - \omega^2_nr_n^2 } \\
    & \lesssim \sum_{j=1}^\infty \frac{2^j t}{2^{2j-2}t^2 - 1 } \hspace{0.2in} [\text{since } r_n^2 \phi_n(r_n^{-1}) \lesssim \sqrt{n}, \ \ \omega_n^2r_n^2 \le 1]\\
    & = \frac{1}{t} \sum_{j=1}^\infty \frac{2^j}{2^{2j-2} - \frac{1}{t^2} } \le \frac{c}{t}
\end{align*}
This proves that: 
$$
r_n \ d\left(\hat \phi, f_0\right) = O_p(1) \,. 
$$
i.e. 
$$
\|\hat f - f_0\|_{L_2(P_X)}^2 = O_p\left(\omega_n^2 + \frac{V_n}{n}\log^{3/2}{n}\right) \,.
$$
Now to balance $\omega_n$ and $V_n$, we choose $N_1 = \lfloor n^{1/2(2\beta_1 + 1)}\rfloor$ and $N_2 = \lfloor n^{1/2(\beta_2 + 1)}\rfloor$. This implies: 
\begin{align*}
\omega_n^2 & \lesssim dn^{-\frac{2\beta_1}{2\beta_1 + 1}} + \binom{d}{2}n^{-\frac{\beta_2}{\beta_2 + 1}} \\
\frac{V_n}{n} & \lesssim \frac{W}{n}\log{W} \lesssim \left(dn^{-\frac{2\beta_1}{2\beta_1 + 1}} + \binom{d}{2}n^{-\frac{\beta_2}{\beta_2 + 1}}\right)\log^3{n},
\end{align*}
which yields
$$
\|\hat f - f_0\|_{L_2(P_X)}^2 = O_p\left(\left(dn^{-\frac{2\beta_1}{2\beta_1 + 1}} + \binom{d}{2}n^{-\frac{\beta_2}{\beta_2 + 1}}\right)\log^{4.5}{n}\right).
$$
This completes the proof. 

\subsection{Proof of Theorem \ref{thm:minimax_lower_bound_growing_d}}
\label{sec:proof_minimax_growing}
Here, we establish the lower bound on the following two-way interaction model: 
$$
    Y_i = \mu + \sum_{j = 1}^d f_j(X_{ij}) + \sum_{k < l} f_{kl}(X_{ik}, X_{il}) + \eps_i \triangleq f_0(X_i) + \eps_i \,.
$$
We show that: 
$$
\inf_{\hat f} \sup_{f \in \cF, X \sim P_{X}} \bbE_{P_{X}}[(\hat f(X) - f(X))^2] \ge c\left(dn^{-\frac{2\beta_1}{2\beta_1 + 1}} + \dbinom{d}{2} n^{-\frac{\beta_2}{\beta_2 + 1}}\right) \,.
$$
where $\cF$ is the collection of all additive functions which satisfies: I) The one-dimensional components are $\beta_1$-smooth, $[0, 1]$ and integrate to $0$ and II) the two-dimensional components are $\beta_2$-smooth, supported on $[0, 1]^2$ having marginals $0$. Our proof is based on the techniques introduced in \cite[Section 2.6.1]{tsybakov2004introduction}. Throughout the proof, assume $X_{ij}\overset{\mathrm{iid}}{\sim} Unif(0,1)$, $i\le n$, $j \le d$ and $\varepsilon_i \overset{\mathrm{iid}}{\sim} N(0,1)$, $i 
\le n$. The key idea is to use Fano's inequality upon carefully choosing a subset of $\cF$. 
Following \cite{tsybakov2004introduction}, fix the following notations: 
\begin{equation}\label{eq:define_mi}
m_i = \lfloor c_0 n^{\frac{1}{2\beta_i + 1}} \rfloor, \quad h_i = m_i^{-1}, \quad i=1,2.    
\end{equation}

\noindent We construct the alternatives in two steps, for one-dimensional components and for two-dimensional components. 
\\\\
{\bf Alternatives for one-dimensional components: } First divide $[0, 1]$ into $m_1$ grids $[(k-1)/m_1, k/m_1]$, $1 \le k \le m_1$. For each $k$, define $x_k = (k - 0.5)/m_1$. Select a smooth kernel $K \in C^\infty(\reals)$ supported on $(-1/2, 1/2)$, $\int K(u) \ du = 0$ and define $\int K^2(u) \ du \triangleq \|K\|_2^2$. Define the functions $\phi_k(x)$, $k \le m_1$ as: 
\begin{equation}\label{eq:define_lower_phi_lin}
    \phi_k(x) := Lh_1^{\beta_1}K\left(\frac{x - x_k}{h_1}\right)
\end{equation}
The constant $L$ is chosen so that $\phi_k$ has bounded derivatives and satisfies Assumption \ref{assn:smoothness}. For each $k$, the function $\phi_k$ is supported on $[(k-1)/m, k/m]$, and the integral of $\phi_k$ is $0$, since
\begin{align*}
    \int \phi_k(x) \ dx & = Lh_1^{\beta_1}\int K\left(\frac{x - x_k}{h_1}\right) \ dx \\
    &=  Lh_1^{\beta_1}\int_{x_k - \frac{h}{2}}^{x_k + \frac{h}{2}} K\left(\frac{x - x_k}{h_1}\right) \ dx \hspace{0.2in} [\text{Since } K(z) > 0 \text{ when } |z| \ge 1/2]\\
    & = Lh_1^{\beta_1 + 1} \int_{-\frac12}^{\frac12} K(z) \ dz = 0 \,.
\end{align*}
For $j \neq k$, $\phi_j$ and $\phi_k$ have disjoint support, yielding the $L_2([0, 1])$ inner product between any $\phi_j$ and $\phi_k$ is $0$. 
Set $M_1 := dm_1$ and $\Omega_1 := \{0, 1\}^{M_1}$. For each $\bomega_1 \in \Omega_1$, define its entries by $\{\bomega_{1, jk}\}_{1 \le j \le d, 1 \le k \le m_1}$. Define functions $f^{(1)}_{\omega_1}$ as: 
$$
f^{(1)}_{\bomega_1}(X) := \sum_{j=1}^d \sum_{k=1}^m \bomega_{1, jk}\phi_k(X_j) \,,
$$
which will contribute to the one dimensional alternatives. 
\\\\
\noindent
{\bf Alternatives for two-dimensional components: } The construction for two-dimensional alternatives are similar to that for one-dimensional components. As before, divide $[0, 1]^2$ into $m_2^2$ rectangles $[(k-1)/m_2, k/m_2] \times [(l-1)/m_2, l/m_2]$, $1 \le k, l \le m_2$. 
Set $x_k := (k - 0.5)/m_2$. Using $K$ same as before, define $\phi_{kl}$, $1 \le k,l \le m_2$ as 
\begin{equation}\label{eq:define_lower_phi_quad}
    \phi_{kl}(x, y) := Lh_2^{\beta_2}K\left(\frac{x - x_k}{h}\right)K\left(\frac{y - x_l}{h}\right)
\end{equation}
Observe that, each $\phi_{kl}$ is supported on $[(k - 1)/m_2, k/m_2] \times [(l-1)/m_2, l/m_2]$. and consequently $L_2([0,1]^2)$ inner-product of any two $\phi_{k_1, l_1}$ and $\phi_{k_2, l_2}$ is 0 if $(k_1,l_1)\neq (k_2,l_2)$. Also the marginals of $\phi_{kl}$ are $0$ as $\int K((x - x_k)/h) \ dx = 0$. 
Set $M_2 := d(d-1)m^2_2/2$ and $\Omega_2 := \{0, 1\}^{M_2}$. For each $\bomega_2 \in \Omega_2$ we write it in a tensor form $\bomega_{2, ijkl}$ where $1 \le i < j \le d$ and $1 \le k \le l \le m_2$ and define $f^{(2)}_{\bomega_2}$ as: 
$$
    f^{(2)}_{\bomega_2}(X)= \sum_{i < j} \sum_{k, l} \bomega_{2, ijkl}\phi_{kl} (X_i, X_j)\,.
$$
These functions will be the two-dimensional component of our alternatives.
\\\\
\noindent
{\bf Final step: } We now construct our alternatives by combining the one-dimensional and two-dimensional components as constructed above. We first set the true mean function $f_0 = 0$. Now we choose $S_1 \subset \Omega_1$ and $S_2 \subset \Omega_2$ carefully and then construct our alternatives as: 
$$
\cF_S = \left\{f_\bomega = f^{(1)}_{\bomega_1} + f^{(2)}_{\bomega_2}, \ \bomega_1 \in S_1, \bomega_2 \in S_2 \right\}
$$
To choose $S_1$ and $S_2$ we invoke Varshamov-Gilbert theorem, cf.~\cite[Lemma 2.9]{tsybakov2004introduction} which we state here for the convenience of the reader:  
\begin{proposition}[Varshamov-Gilbert]
\label{prop:VG}
    For any $M\ge 8$, there exists $S \subseteq \{0, 1\}^M$ with $|S| \ge 2^{M/8}$, such that for any $\bomega \neq \bomega' \in S$, $\rho(\bomega, \bomega') \ge M/8$, where $\rho$ is the Hamming distance. 
\end{proposition}
\noindent By Proposition \ref{prop:VG}, we can choose $S_1$ and $S_2$ such that for any $\bomega_1, \bomega_1' \in S_1$ we have $\rho(\bomega_1, \bomega_1') \ge M_1/8$ and for any $\bomega_2, \bomega_2' \in S_2$, $\rho(\bomega_2, \bomega_2') \ge M_2/8$. Observe that for any $\bomega = (\bomega_1, \bomega_2) \neq \bomega' = (\bomega'_1, \bomega_2')$ we have:  
\begin{align*}
    & d^2\left(f_\bomega, f_{\bomega'}\right) = \int_{[0, 1]^d} (f_\bomega(\bx) - f_{\bomega'}(\bx))^2 \ d\bx  \\
    & = \int_{[0, 1]^d} \left(\sum_{j=1}^d \sum_{k=1}^m (\bomega_{1, jk} - \bomega'_{1, jk})\phi_k(X_j) + \sum_{i < j} \sum_{k, l} (\bomega_{2, ijkl} - \bomega'_{2, ijkl})\phi_{kl}(X_i, X_j)\right)^2 \ dx \\
    & = \int_{[0, 1]^d} \left(\sum_{j=1}^d \sum_{k=1}^m (\bomega_{1, jk} - \bomega'_{1, jk})\phi_k(X_j) \right)^2 \ d\bx + \left(\sum_{i < j} \sum_{k, l} (\bomega_{2, ijkl} - \bomega'_{2, ijkl})\phi_{kl}(X_i, X_j)\right)^2 \ d\bx \\
    & \qquad \qquad + 2\int_{[0, 1]^d} \left(\sum_{j=1}^d \sum_{k=1}^m (\bomega_{1, jk} - \bomega'_{1, jk})\phi_k(X_j) \right)\left(\sum_{i < j} \sum_{k, l} (\bomega_{2, ijkl} - \bomega'_{2, ijkl})\phi_{kl}(X_i, X_j)\right) \ d\bx \\
    & \triangleq T_1 + T_2 + 2T_3 \,.
\end{align*}
We now analyze each $T_i$ separately.
\begin{align*}
    T_1 &= \int_{[0, 1]^d} \left(\sum_{j=1}^d \sum_{k=1}^m (\bomega_{1, jk} - \bomega'_{1, jk})\phi_k(X_j) \right)^2 \ d\bx \\
    & = \sum_{j=1}^d \sum_{k=1}^m (\bomega_{1, jk} - \bomega'_{1, jk})^2\int_0^1 \phi^2_k(X_j) \ dX_j \\
    & = L^2 h_1^{2\beta_1}\sum_{j=1}^d \sum_{k=1}^m (\bomega_{1, jk} - \bomega'_{1, jk})^2\int_0^1 K^2\left(\frac{X_j - x_k}{h_1}\right) \ dX_j \\
    & = L^2 h_1^{2\beta_1 + 1}\sum_{j=1}^d \sum_{k=1}^m (\bomega_{1, jk} - \bomega'_{1, jk})^2\int_{-\frac12}^{\frac12} K^2\left(z\right) \ dz \\
    & = L^2 h_1^{2\beta_1 + 1}\|K\|_2^2 \ \sum_{j=1}^d \sum_{k=1}^m (\bomega_{1, jk} - \bomega'_{1, jk})^2 = L^2 h_1^{2\beta_1 + 1}\|K\|_2^2 \ \rho(\bomega_1, \bomega'_1) 
\end{align*}
For $T_2$: 
\begin{align*}
    T_2 & = \int_{[0, 1]^d}\left(\sum_{i < j} \sum_{k, l} (\bomega_{2, ijkl} - \bomega'_{2, ijkl})\phi_{kl}(X_i, X_j)\right)^2 \ d\bx \\
    & = \sum_{i < j} \sum_{k, l} (\bomega_{2, ijkl} - \bomega'_{2, ijkl})^2 \int_{[0, 1]^2} \phi^2_{kl}(X_i, X_j) \ dX_i \ dX_j \\
    & = L^2h_2^{2\beta_2}\sum_{i < j} \sum_{k, l} (\bomega_{2, ijkl} - \bomega'_{2, ijkl})^2 \int_{[0, 1]^2} K^2\left(\frac{X_i - x_j}{h_2}\right)K^2\left(\frac{X_j - x_k}{h_2}\right) \ dX_i \ dX_j \\
    & =  L^2h_2^{2\beta_2 + 2} \sum_{i < j}\sum_{k, l} (\bomega_{2, ijkl} - \bomega'_{2, ijkl})^2 \left(\int_{-\frac12}^{\frac12} K^2(z) \ dz\right)^2 \\
    & = L^2h_2^{2\beta_2 + 2}\|K\|_2^4 \ \sum_{i < j}\sum_{k, l} (\bomega_{2, ijkl} - \bomega'_{2, ijkl})^2 = L^2h_2^{4\beta_2 + 2}\|K\|_2^4 \rho(\bomega_2, \bomega_2') \,.
\end{align*}
Finally we claim that $T_3 = 0$. This follows from the fact that if $i, j, k$ all are distinct, then for any $1 \le v_1 \le m_1$ and $1 \le v_2, v_3 \le m_2$: 
\begin{align*}
    & \int_{[0, 1]^3}\phi_{v_1}(X_i) \phi_{v_2, v_3}(X_j, X_k) \ dX_i \ dX_j \ dX_k \\
    & = \int_{[0, 1]}\phi_{v_1}(X_i) \ dX_i \ \int_{[0, 1]^2}\phi_{v_2, v_3}(X_j, X_k) \ dX_j \ dX_k = 0
\end{align*}
as the integrals of $\phi_{v_1}$ and $\phi_{v_2, v_3}$ are $0$ by construction. When they are not all distinct, suppose $i =j \neq k$. Then: 
\begin{align*}
    & \int_{[0, 1]^2}\phi_{v_1}(X_i) \phi_{v_2, v_3}(X_i, X_k) \ dX_i \ dX_k \\
    & = \int_{[0, 1]}\phi_{v_1}(X_i) \ \left(\int_{[0, 1]}\phi_{v_2, v_3}(X_i, X_k) \ dX_k\right) \ dX_i = 0
\end{align*}
as marginals of $\phi_{v_2, v_3}$ is 0 by our construction. Therefore $T_3 = 0$. Combining the expression of $T_i$'s, we have: 
\begin{align}
\label{eq:d_val}
d^2(f_\bomega, f_{\bomega'}) & =  L^2 h_1^{2\beta_1 + 1}\|K\|_2^4 \ \rho(\bomega_1, \bomega_1') +  L^2h_2^{2\beta_2 + 2} \|K\|_2^4\ \rho(\bomega_2, \bomega_2')
\end{align}
where $\rho$ is the Hamming distance. 
From Proposition \ref{prop:VG}, we have: 
\begin{align}
\label{eq:d_val_lb}
d^2(f_\bomega, f_{\bomega'}) & \ge  L_1^2 h_1^{2\beta_1 + 1}\|K\|_2^2 \ \frac{dm_1}{8} +  L_2^2h_2^{2\beta_2 + 2} \|K\|_2^4\ \frac{d(d-1)m_2^2}{16} \notag \\
& =  \frac{L_1^2 \|K\|_2^2}{8} dm_1^{-2\beta_1} +  \frac{L_2^2 \|K\|_2^4}{8}\dbinom{d}{2}m_2^{-2\beta_2} := 4\delta^2 
 \,.
\end{align}
Here the second equality follows from the fact that $h_1 = m_1^{-1}$ and $h_2= m_2^{-1}$ and the third inequality follows from the definition of $m_1, m_2$. 
Applying Fano's inequality~\cite[Proof of Theorem 2.18]{mukherjee2021optimal} on the collection $\cF_S$, we obtain: 
\begin{equation}
    \label{eq:Fano}
    \inf_{\hat f} \sup_{f,P_{X}} \bbE_{P_{X}}[(\hat f(X) - f(X))^2] \ge \delta^2 \left(1 -\frac{\frac{n}{|S|^2} \sum_{\omega \neq \omega' \in S}KL\left(\bbP_\omega | \bbP_{\omega'}\right) + \log{2}}{\log{(|S| - 1)}}\right)
\end{equation}
As the error $\eps_i$'s are normal, we have: 
\begin{align*}
    KL\left(\bbP_\omega | \bbP_{\omega'}\right) & = \frac{1}{2}\bbE\left[\left(f_\omega(X) - f_{\omega'}(X)\right)^2\right] = \frac12 d^2\left(f_\bomega, f_{\bomega'}\right) \\ 
    & =  L^2 h_1^{2\beta_1 + 1}\|K\|_2^2 \ \rho(\bomega_1, \bomega_1') +  L^2h_2^{2\beta_2 + 2} \|K\|_2^4\ \rho(\bomega_2, \bomega_2') \\
    & \le L^2 h_1^{2\beta_1 + 1}\|K\|_2^2 \ M_1+  L^2h_2^{2\beta_2 + 2} \|K\|_2^4\ M_2 \\
    & = L^2 h_1^{2\beta_1 + 1}\|K\|_2^2 \ dm_1+  L^2h_2^{2\beta_2 + 2} \|K\|_2^4\ \dbinom{d}{2} m_2^2 \\
    & = L^2\|K\|_2^2 dm_1^{-2\beta_1}+  L^2\|K\|_2^4\ \dbinom{d}{2} m_2^{-2\beta_2} \le 32 \delta^2 \,.
\end{align*}
Moreover, $|S| = |S_1 \times S_2| \ge 2^{(M_1 + M_2)/8}$ by Proposition \ref{prop:VG}. Plugging the bounds in \eqref{eq:Fano}, we have: 
\begin{align}
\label{eq:FANO_2}
    \inf_{\hat f} \sup_{f,P_{X}} \bbE_{P_{X}}[(\hat f(X) - f(X))^2] \ge \delta^2 \left(1 -9\frac{32n\delta^2 + \log{2}}{M_1 + M_2}\right)
\end{align}
Now from the definition of $m_1, m_2$ we have: 
\begin{align*}
    n\delta^2 & = n\left(\frac{L_1^2 \|K\|_2^2}{32} dm_1^{-2\beta_1} +  \frac{L_2^2 \|K\|_2^4}{32}\dbinom{d}{2}m_2^{-2\beta_2}\right) \\
    & \le Cn \left(\frac{L_1^2 \|K\|_2^2}{32} dn^{-\frac{2\beta_1}{2\beta_1 + 1}} +  \frac{L_2^2 \|K\|_2^4}{32}\dbinom{d}{2}n^{-\frac{2\beta_2}{2\beta_2 + 2}}\right)  \\
    & = C\left(\frac{L_1^2 \|K\|_2^2}{32} \vee \frac{L_2^2 \|K\|_2^4}{32}\right)\left( dn^{\frac{1}{2\beta_1 + 1}} +  \dbinom{d}{2}n^{\frac{2}{2\beta_2 + 2}}\right) \\
    & \triangleq C\left(\frac{L_1^2 \|K\|_2^2}{32} \vee \frac{L_2^2 \|K\|_2^4}{32}\right)\psi_n \,.
\end{align*}
On the other hand: 
\begin{align*}
    M_1 + M_2 & = dm_1 + \dbinom{d}{2}m_2^2 \ge c\left(dn^{\frac{1}{2\beta_1 + 1}} +  \dbinom{d}{2}n^{\frac{2}{2\beta_2 + 2}}\right) = c\psi_n \,.
\end{align*}
Using these bounds in \eqref{eq:FANO_2} yields: 
\begin{align*}
    \inf_{\hat f} \sup_{f,P_{X}} \bbE_{P_{X}}[(\hat f(X) - f(X))^2] & \ge \delta^2 \left(1 -9\frac{C\left(\frac{L_1^2 \|K\|_2^2}{32} \vee \frac{L_2^2 \|K\|_2^4}{32}\right)\psi_n + \log{2}}{c\psi_n}\right) \\
    & = \delta^2 \left(1 -9\frac{C\left(\frac{L_1^2 \|K\|_2^2}{32} \vee \frac{L_2^2 \|K\|_2^4}{32}\right) + \frac{\log{2}}{\psi_n}}{c}\right)
\end{align*}
The constants $C, c$ depends on $c_1, c_2$ via the definition of $m_1,m_2$. Choosing them appropriately and using the fact $\psi_n \uparrow \infty$ as $n \uparrow \infty$, we can make 
$$
9\frac{C\left(\frac{L_1^2 \|K\|_2^2}{32} \vee \frac{L_2^2 \|K\|_2^4}{32}\right) + \frac{\log{2}}{\psi_n}}{c} \le \frac12 \,.
$$
which implies: 
\begin{align*}
 \inf_{\hat f} \sup_{f,P_{X}} \bbE_{P_{X}}[(\hat f(X) - f(X))^2] & \ge \delta^2 \\
 & \ge \left(\frac{L_1^2 \|K\|_2^2}{8} \wedge \frac{L_2^2 \|K\|_2^4}{8}\right) \left( dn^{\frac{-2\beta_1}{2\beta_1 + 1}} +  \dbinom{d}{2}n^{\frac{-2\beta_2}{2\beta_2 + 2}}\right) \,.
\end{align*}
This completes the proof.

\subsection{Proof of Theorem \ref{thm:fixed_design_wo_rsc} and \ref{thm:fixed_error_sparse}}\label{sec:proof_fixed_high_dim}
Recall the for any $\phi = \sum_j \phi_j + \sum_{k < l} \phi_{kl}$, we define $\|\phi\|_{n, 1} := \sum_j \|\phi_j\|_{n} + \sum_{j < k}\|\phi_{jk}\|_{n}$ and $\|\phi\|^2_n =\sum_{i=1}^{n} \phi^2(X_i)/n$. Further define $\|\phi\|_{n,lin}:=\sum_j \|\phi_j\|_{n}$ and $\|\phi\|_{n,quad}:=\sum_{k < l} \|\phi_{kl}\|_{n}$. From Theorem \ref{thm:statistical_error_lowdim}, there exists $\{\phi^\star_j\}_{j \in S_1}, \left\{\phi^\star_{kl}\right\}_{(j, k) \in S_2}$, $j\in S_1, (j,k) \in S_2$, such that, $\phi^\star_j \in \mathcal{F}_{NN,1}$, $\phi^\star_{jk} \in \mathcal{F}_{NN,2}$ and
\begin{align}
    & \|f_{0,j} - \phi^\star_{j}\|_\infty \le C_1 (N_1L_1)^{-2\beta_1} \ \forall \ j \in S_1 \,, \label{eq:nn_approx_lin} \\
    & \|f_{0,kl} - \phi^\star_{kl}\|_\infty \le C_2 (N_2L_2)^{-\beta_2} \ \forall \ (k < l) \in S_2 \,.\label{eq:nn_approx_quad}
\end{align}
and furthermore, $\phi^\star_j$ has integral $0$ and $\phi^\star_{kl}$ has $0$ marginals. For the simplicity of the rest of the proof, take $L_i$ to be of the constant order. Therefore the architecture we use here is the following: for the fitting the univariate function we use DNNs from the class $\cF(1, c_1 N_1 \log{N_1}, c_2, 1)$ and for fitting bivariate functions, we use DNNs from the class  $\cF(1, c_1 N_2 \log{N_2}, c_3, 1)$. Therefore, defining $\phi^\star = \sum_{j \in S_1} \phi^\star_j + \sum_{(j, k) \in S_2} \phi^\star_{jk}$, we have by Theorem \ref{thm:approx_growing_d}:  
\begin{align}
    \bbE\left[(f_0(X)-\phi^\star(X))^2\right] & \le C_1 \left(s_1 (N_1)^{-4\beta_1}+ s_2 (N_2)^{- 2\beta_2}\right) \notag \\
    & \triangleq C_1(s_1 \rho_{n, 1}^2+ s_2 \rho_{n, 2}^2)\,.\label{eq:define_nn_approx_linf}
\end{align}



\noindent As $\hat \phi$ is the minimizer of the penalized loss function among the class of neural networks, it outperforms $\phi^\star$. Setting $pen(\phi):= \lambda_{n,1} \| \phi\|_{n, lin} +\lambda_{n,2} \|\phi\|_{n, quad}$, we have: 
\begin{align}
    & \frac{1}{2}\|Y - \hat \phi\|_n^2 + pen(\hat\phi) \le  \frac{1}{2}\|Y - \phi^\star\|_n^2 + pen(\phi^\star) \notag \\
    \implies &  \frac{1}{2}\|\hat \phi - \phi^\star\|_n^2 + pen(\hat\phi) \le  \langle f_0 - \phi^\star, \hat \phi - \phi^\star\rangle_n + \langle \eps , \hat \phi - \phi^\star\rangle_n + pen(\phi^\star) \notag \\
    \implies & \frac{1}{2}\|\hat \phi - \phi^\star\|_n^2 + pen(\hat\phi)\le \|f_0 - \phi^\star\|_n^2 + \frac14 \|\hat \phi - \phi^\star\|_n^2 + \langle \eps , \hat \phi - \phi^\star\rangle_n + pen(\phi^\star) \notag \\
    \implies & \frac{1}{4}\|\hat \phi - \phi^\star\|_n^2 + pen(\hat\phi) \le \|f_0 - \phi^\star\|_n^2 + \langle \eps , \hat \phi - \phi^\star\rangle_n +pen(\phi^\star). \label{eq:fan1}
\end{align}
\noindent
First we bound the empirical error $\|f_0 - \phi^\star\|_n^2$ in terms of its population counterpart $\|f_0 - \phi^\star\|_2^2$. Using Chebychev inequality, we have for $t > 0$: 
\begin{equation}
    \label{eq:chebychev}
    \bbP\left(\left|\|f_0 - \phi^\star\|_n^2 - \|f_0 - \phi^\star\|_2^2  \right| > t\right) \le \frac{\bbE\left[(f_0 - \phi^\star)^4\right]}{nt^2}
\end{equation}
Expanding the fourth moment, we obtain: 
\begin{align*}
    & \bbE\left[(f_0 - \phi^\star)^4\right] \le 8\left\{\bbE\left[\left(\sum_{j \in S_1} (f_{0,j} - \phi^\star_j)\right)^4\right] + \bbE\left[\left(\sum_{(j, k) \in S_2} (f_{0,jk} - \phi^\star_{jk})\right)^4\right]\right\} \\
    & \le 8p_{\max}\left\{\int_{[0, 1]^d}\left(\sum_{j \in S_1} (f_{0,j}(x_j) - \phi^\star_j(x_j))\right)^4 \ d\bx + \int_{[0, 1]^d}\left(\sum_{(j, k) \in S_2} (f_{0,jk}(x_j, x_k) - \phi^\star_{jk}(x_j, x_k))\right)^4 \ d\bx \right\} \\
    & \triangleq 8p_{\max}(T_1 + T_2) \,.
\end{align*}
Now we analyze $T_1$ and $T_2$ separately. For $T_1$: 
\begin{align*}
    T_1 & = \int_{[0, 1]^d}\left(\sum_{j \in S_1} (f_{0,j}(x_j) - \phi^\star_j(x_j))\right)^4 \ d\bx \\
    & = \sum_{j \in S_1} \int_0^1 (f_{0,j}(x_j) - \phi^\star_j(x_j))^4 \ dx_j + \sum_{j \neq k \in S_1} \int_0^1 (f_{0,j}(x_j) - \phi^\star_j(x_j))^2 \ dx_j\int_0^1 (f_{0,j}(x_k) - \phi^\star_k(x_k))^2 \ dx_k \\
    & \le C^4 s^2_1 (N_1L_1)^{-8\beta} \,.
\end{align*}
where the last line follows from the fact that $\|f_{0, j} - \phi^\star_j\|_\infty \le C(N_1L_1)^{-2\beta}$. 
Similar analysis for $T_2$ yields: 
\begin{align*}
     T_2 & = \int_{[0, 1]^d}\left(\sum_{(j, k) \in S_2} (f_{0,jk}(x_j, x_k) - \phi^\star_{jk}(x_j, x_k))\right)^4 \ d\bx  
    \le C^4 s^2_2 (N_2L_2)^{-4\beta} \,.
\end{align*}
using $\|f_{0, jk} - \phi^\star_{jk}\|_\infty \le C(N_2L_2)^{-\beta}$,  repeatedly and using Assumption \ref{asm:d_growth_high}. Using the bounds of $T_1$ of $T_2$, we have:  
\begin{align*}
    \bbE\left[(f_0 - \phi^\star)^4\right] &\le 8C^4p_{\max}\left(s^2_1 (N_1L_1)^{-8\beta} + s^2_2 (N_2L_2)^{-4\beta}\right)\\
    &\le 8C^4p_{\max} \Big(s^2_1\rho^4_{n,1}+s^2_2\rho^4_{n,2} \Big)\,.
\end{align*}
Hence, choosing $t_0 = C^2 \sqrt{\frac{8p_{\max} (s^2_1\rho^4_{n,1}+s^2_2\rho^4_{n,2}) \log{n}}{n}}$, we have from \eqref{eq:chebychev}: 
\begin{equation}
    \label{eq:chebychev_2}
    \bbP\left(\left|\|f_0 - \phi^\star\|_n^2 - \|f_0 - \phi^\star\|_2^2  \right| > t_0\right) \le \frac{1}{\log{n}} 
\end{equation}
Define the event $\Omega_{n, 1}$ as: 
$$
\Omega_{n, 1}  = \left\{\left|\|f_0 - \phi^\star\|_n^2 - \|f_0 - \phi^\star\|_2^2  \right| \le t_0 \right\} \,.
$$
Noticing $t_0 \ll s_1 \rho_{n, 1}^2+ s_2 \rho_{n, 2}^2$, by \eqref{eq:define_nn_approx_linf}, we have 
on $\Omega_{n, 1}$ that
\begin{align*}
    \eqref{eq:fan1}  \implies &  \frac{1}{4}\|\hat \phi - \phi^\star\|_n^2 + pen(\hat \phi) \le C_3(s_1 \rho_{n, 1}^2+ s_2 \rho_{n, 2}^2)  + \langle \eps , \hat \phi - \phi^\star\rangle_n + pen(\phi^\star)
\end{align*}
In the next step, we bound the empirical error, i.e. the inner product between $\eps$ and $\hat \phi - \phi^\star$. First of all, note that, by triangle inequality: 
\begin{align*}
\left|\langle \eps , \hat \phi - \phi^\star\rangle_n \right| & \le 
\sum_{j=1}^d \left|\frac1n \sum_i \eps_i \left(\hat \phi_j(X_{ij}) - \phi^\star_j(X_{ij})\right)\right| \\
& \qquad \qquad + \sum_{k<l} \left|\frac1n \sum_i \eps_i \left(\hat \phi_{kl}(X_{ik}, X_{il}) - \phi^\star_{kl}(X_{ik}, X_{il})\right)\right|
\end{align*}
To bound each of the summands on the RHS, we first bound each term inside the summation. Towards that direction, we use Lemma $4$ of \cite{fan2022factor}, which says that if $\cG_n$ is a class of functions with VC dimension $V_n$, then for any fixed $g_0 \in\cG, \eps> 0, t > 0$, with probability $\ge 1 - \log{(1/\eps)}e^{-t}$: 
\begin{align*}
\frac{1}{n}\left|\sum_{i=1}^n\eps_i(g(X_i) - g_0(X_i))\right| & \le C_4\left(\|g - g_0\|_n + \eps\right)\sqrt{\frac{V_{n, 1} \log{n}}{n} + \frac{t}{n}} 
\end{align*}
for some universal constant $c$. Note that, $\phi_j \in \cF_{NN}(1, c_1 N_1 \log{N_1}, c_2, c_3 (N_1 \log{N_1})^2, 1)$, whose VC dim $V_{n, 1} \le C_5 N_1^2\log{N_1})$ (see Lemma \ref{lem:bartlett_vc_bound}). 
We now apply this lemma to each of the component functions with $t = 2\log{d}$ and use a union bound to conclude:
\begin{align}\label{eq:good_event_n2_1}
  \sum_{j=1}^d \left|\frac1n \sum_i \eps_i \left(\hat \phi_j(X_{ij}) - \phi^\star_j(X_{ij})\right)\right| \le  C_4\sqrt{\frac{V_{n, 1} \log{p}}{n} + \frac{2\log{d}}{n}} \left(\sum_{j=1}^d\|\hat \phi_j - \phi^\star_{n, j}\|_n + d\eps_1\right)\,.
\end{align}
The above event occurs with probability $\ge 1 - \log(1/\eps_1)e^{\log{d}-2\log{d}} = 1 - \log(1/\eps_1)/d$. A similar calculation for the bivariate components yields: 
\begin{align}
\label{eq:good_event_n2_2}
    & \sum_{k<l} \left|\frac1n \sum_i \eps_i \left(\hat \phi_{kl}(X_{ik}, X_{il}) - \phi^\star_{kl}(X_{ik}, X_{il})\right)\right| \notag \\
    & \hspace{5em}\le C_5\sqrt{\frac{V_{n, 2} \log{n}}{n} + \frac{3\log{d}}{n}} \left(\sum_{k < l}\|\hat \phi_{kl} - \phi^\star_{kl}\|_n + \frac{d(d-1)}{2}\eps_2\right) \,,
\end{align}
with probability $\ge 1 - \log{(1/\eps_2)}e^{2\log{d}-3\log{d}} = 1 - \log{(1/\eps_2)}/d$. For the rest of the analysis, define the penalty parameters $\lambda_{n, 1}, \lambda_{n, 2}$ and $\eps_1, \eps_2$ as: 
\begin{align}
\label{eq:lambda_1_val} \lambda_{n, 1} & = 2C_4\sqrt{\frac{V_{n, 1} \log{p}}{n} + \frac{2\log{d}}{n}}\\
\label{eq:lambda_2_val} \lambda_{n, 2} & = 2C_5\sqrt{\frac{V_{n, 2} \log{n}}{n} + \frac{3\log{d}}{n}}  \\
\label{eq:eps_val} \eps_i & = \frac{s_i \lambda_{n, i}}{\binom{d}{i}} \text{ for } i=1,2.
\end{align}
Combining equation \eqref{eq:good_event_n2_1} and \eqref{eq:good_event_n2_2} we have: 
\begin{align}
    \label{eq:good_event_n2_3}
    \left|\langle \eps, \hat \phi - \phi^\star \rangle_n\right| & \le \frac{\lambda_{n, 1}}{2} \sum_{j=1}^p\|\hat \phi_j - \phi^\star_{n, j}\|_{n, lin} +  \frac{\lambda_{n, 2}}{2} \sum_{k < l}\|\hat \phi_{kl} - \phi^\star_{kl}\|_{n, quad} \notag \\
    & \hspace{15em} + \frac{s_1\lambda_{n, 1}^2}{2} + \frac{s_2\lambda_{n, 2}^2}{2} \,.
\end{align}
Define this event \eqref{eq:good_event_n2_3} to be $\Omega_{n, 2}$ and we will later show that our choice of $(N_1, N_2)$ ensures that $\bbP(\Omega_{n, 2}) \to 1$. On the event $\Omega_{n, 1} \cap \Omega_{n, 2}$, we have: 
\begin{align}
& \frac{1}{4}\|\hat \phi - \phi^\star\|_n^2 + \lambda_{n, 1} \|\hat \phi \|_{n, lin} + \lambda_{n, 2}\|\hat \phi\|_{n, quad} \notag \\
& \le C_3 (s_1\rho_{n, 1}^2 + s_2 \rho_{n, 2}^2)  + \frac{\lambda_{n, 1}}{2} \sum_{j=1}^p\|\hat \phi_j - \phi^\star_{n, j}\|_{n, lin} \notag \\
& \qquad \qquad +  \frac{\lambda_{n, 2}}{2} \sum_{k < l}\|\hat \phi_{kl} - \phi^\star_{kl}\|_{n, quad}  + \frac{s_1\lambda_{n, 1}^2}{2} + \frac{s_2\lambda_{n, 2}^2}{2} + \lambda_{n, 1} \|\hat \phi \|_{n, lin} + \lambda_{n, 2}\|\hat \phi\|_{n, quad} \notag 
\end{align}
Some simple algebra yields: 
\begin{align}
\label{eq:final_bound_fixed}
    & \frac{1}{4}\|\hat \phi - \phi^\star\|_n^2 + \frac{\lambda_{n, 1}}{2} \|\hat \phi \|_{n, lin} + \frac{\lambda_{n, 2}}{2}\|\hat \phi\|_{n, quad} \notag \\
    & \le C_3 (s_1\rho_{n, 1}^2 + s_2 \rho_{n, 2}^2)  + \frac{3\lambda_{n, 1}}{2} \sum_{j \in S_1}\|\hat \phi_j - \phi^\star_{n, j}\|_{n, lin} +  \frac{3\lambda_{n, 2}}{2} \sum_{k < l \in S_2}\|\hat \phi_{kl} - \phi^\star_{kl}\|_{n, quad} \notag \\
    & \hspace{5em} + \frac{s_1\lambda_{n, 1}^2}{2} + \frac{s_2\lambda_{n, 2}^2}{2} \,.
\end{align}
If we don't assume any \emph{Restricted Strong Convexity}(RSC) on the underlying function class (i.e. in our case DNNs with some pre-specified architecture), then we can use the fact $\|\hat \phi_j - \phi^\star_{n, j}\|_n \le 2B$ for $j$ and consequently we have: 
\begin{align}
    \label{eq:bound_fixed_no_rse}
    & \frac{1}{4}\|\hat \phi - \phi^\star\|_n^2 + \frac{\lambda_{n, 1}}{2} \|\hat \phi \|_{n, lin} + \frac{\lambda_{n, 2}}{2}\|\hat \phi\|_{n, quad} \notag \\
    & \le C_3 (s_1\rho_{n, 1}^2 + s_2 \rho_{n, 2}^2)  + \frac{3}{2}(s_1\lambda_{n, 1} + s_2\lambda_{n, 2}) + \frac{s_1\lambda_{n, 1}^2}{2} + \frac{s_2\lambda_{n, 2}^2}{2} \,.
\end{align}
This completes the proof of Theorem \ref{thm:fixed_design_wo_rsc}. For the second part, we assume to have RSC condition 
 Assumption \eqref{assn:rsc}. From \eqref{eq:final_bound_fixed}, we obtain: 
\begin{align*}
    & \frac{1}{4}\|\hat \phi - \phi^\star\|_n^2 + \frac{\lambda_{n, 1}}{2} \|\hat \phi \|_{n, lin} + \frac{\lambda_{n, 2}}{2}\|\hat \phi\|_{n, quad} \notag \\
    & \le C_3 (s_1\rho_{n, 1}^2 + s_2 \rho_{n, 2}^2)  + \frac{3\lambda_{n, 1}}{2} \sum_{j \in S_1}\|\hat \phi_j - \phi^\star_{n, j}\|_{n, lin} +  \frac{3\lambda_{n, 2}}{2} \sum_{k < l \in S_2}\|\hat \phi_{kl} - \phi^\star_{kl}\|_{n, quad} \notag \\
    & \hspace{5em} + \frac{s_1\lambda_{n, 1}^2}{2} + \frac{s_2\lambda_{n, 2}^2}{2} \\
    & \le C_3(s_1(\rho_{n, 1}^2 + \lambda_{n, 1}^2) + s_2(\rho_{n, 2}^2 + \lambda_{n, 2}^2)) + \frac{3}{2}\lambda_{n,1}\sqrt{s_1}\sqrt{\sum_{j \in S_1}\|\hat \phi_j - \phi^\star_{n, j}\|^2_{n, lin}} \\
    & \qquad \qquad \qquad + \frac{3}{2} \lambda_{n, 2}\sqrt{s_2}\sqrt{\sum_{k < l \in S_2}\|\hat \phi_{kl} - \phi^\star_{kl}\|^2_{n, quad}} \\
    & \le C_3 (s_1(\rho_{n, 1}^2 + \lambda_{n, 1}^2) + s_2(\rho_{n, 2}^2 + \lambda_{n, 2}^2)) + \frac{9s_1\lambda_{n, 1}^2}{2\kappa_1^2} + \frac{\kappa_1^2}{8}\sum_{j \in S_1}\|\hat \phi_j - \phi^\star_{n, j}\|^2_{n, lin} \\
    & \qquad \qquad \qquad + \frac{9s_2\lambda_{n, 2}^2}{2\kappa_2^2} + \frac{\kappa_2^2}{8}\sum_{k < l \in S_2}\|\hat \phi_{kl} - \phi^\star_{kl}\|^2_{n, quad} \\
    & \le C_3(s_1(\rho_{n, 1}^2 + \lambda_{n, 1}^2) + s_2(\rho_{n, 2}^2 + \lambda_{n, 2}^2)) +  \frac{9s_2\lambda_{n, 2}^2}{2\kappa_2^2} + \frac{9s_1\lambda_{n, 1}^2}{2\kappa_1^2} + \frac{1}{8}\|\hat \phi - \phi^\star\|_n^2 
\end{align*}
This implies: 
\begin{align*}
& \frac{1}{8}\|\hat \phi - \phi^\star\|_n^2 + \frac{\lambda_{n, 1}}{2} \sum_{j \in S_1^c}\|\hat \phi_j \|_{n} + \frac{\lambda_{n, 2}}{2}\sum_{(j < k) \in S_2^c}\|\hat \phi_{jk}\|_{n} 
& \le C_6\left(s_1(\rho_{n, 1}^2 + \lambda_{n, 1}^2) + s_2(\rho_{n, 2}^2 + \lambda_{n, 2}^2)\right)
\end{align*}
Now we select the optimal value of $N_i$'s by balancing $\rho_{n, i}$ and $\lambda_{n, i}$. First take $i = 1$. Recall that we have: 
$$
\rho_{n, 1}^2 = N_1^{-4\beta_1}, \text{   and   } \lambda^2_{n, 1} = 4C_4^2 \left(\frac{V_{n, 1}\log{n}}{n} + \frac{2\log{d}}{n}\right) \le C_7 \left(\frac{N_1^2\log^3{N_1}\log{n}}{n}  + \frac{\log{d}}{n}\right)
$$
Choosing $N_1 = \lfloor n^{\frac{1}{2(2\beta_1 + 1)}}\rfloor$ yields: 
\begin{align}
\label{eq:rho_val_final}
    \rho^2_{n, 1} & \sim n^{-\frac{2\beta_1}{2\beta_1 + 1}} \,, \\
\label{eq:lambda_val_final}
    \lambda^2_{n, 1} & \sim n^{-\frac{2\beta_1}{2\beta_1 + 1}}\log^4{n} + \frac{\log{d}}{n} \,.
\end{align}
In particular we have: 
$$
\rho_{n, 1}^2 + \lambda^2_{n, 1} \le C_8 \left(n^{-\frac{2\beta_1}{2\beta_1 + 1}}\log^4{n} + \frac{\log{d}}{n}\right) \,.
$$
Similar calculation for $i = 2$ implies that choosing $N_2 = \lfloor n^{\frac{1}{2(\beta_2 + 1)}}\rfloor$ we have: 
$$
\rho_{n, 2}^2 + \lambda^2_{n, 2} \le C_9 \left(n^{-\frac{\beta_2}{\beta_2 + 1}}\log^4{n} + \frac{\log{d}}{n}\right) \,.
$$
Therefore, the above choice of penalty parameters yields: 
\begin{align*}
& \frac{1}{8}\|\hat \phi - \phi^\star\|_n^2 + \frac{\lambda_{n, 1}}{2} \sum_{j \in S_1^c}\|\hat \phi_j \|_{n} + \frac{\lambda_{n, 2}}{2}\sum_{(j < k) \in S_2^c}\|\hat \phi_{jk}\|_{n} \\
& \qquad \qquad \le C_{10}\left(s_1\left(n^{-\frac{2\beta_1}{2\beta_1 + 1}}\log^4{n} + \frac{\log{d}}{n}\right) + s_2\left(n^{-\frac{\beta_2}{\beta_2 + 1}}\log^4{n} + \frac{\log{d}}{n}\right)\right) \,.
\end{align*}
This completes the proof.

\subsection{Proof of Theorem \ref{thm:random_design_hard_threshold_ubd}}
\label{sec:proof_ht}
For notational simplicty, define $r^2_n = \tilde C(s_1\lambda^2_{n, 1} + s_2\lambda^2_{n, 2})$ for some constant $\tilde C$ mentioned explicitly later in the proof. First we show that the set $\hat S_i$ for $i \in \{1, 2\}$, obtained by hard thresholding (Step 3 of Algorithm \ref{algo:random_design_algo}) satisfies $\hat S_i \supset S_i$ with high probability (whp). Recall that $\|f_j^0\|_2 > r_n$ and $\|f_{jk}^0\|_2 > r_n$ for all $j \in S_1$ and $(j < k) \in S_2$. Using \cite{lu2021deep}, there exists DNNs $\{\phi^\star_j\}_j$ and $\{\phi^\star_{jk}\}$ such that $\|\phi^\star_j\|_2 > r_n/2$ and $\|\phi^\star_{jk}\|_2 > r_n/2$ as the approximation error for each component is $< r_n/2$ by the definition of $r_n$ for all $j \in S_1$ and $(j < k) \in S_2$. We will show first that $S_1 \subseteq \hat S_1$ whp. The key idea is as follows: if $S_1 \nsubseteq \hat S_1$, there exists $j \in S_1$ such that $j \notin \hat S_1$, i.e. $\|\hat \phi_j\|_n < C_2 \lambda_{n, 1}$. Since $\|\phi^\star_j\|_2 \ge r_n/2$, we have $\|\phi^\star_j\|_n \ge r_n/4$ whp (which again follows from the fact that $nr_n^2 \to \infty$, details later). If $\|\hat \phi_j\|_n < C_2 \lambda_{n, 1}$, then 
\begin{equation}
\label{eq:fixed_dim_prop}
C(s_1 \lambda_{n, 1}^2 + s_2 \lambda_{n, 2}^2) \ge \frac{1}{8}\|\hat \phi - \phi^\star\|^2_n \ge \frac18 \|\hat \phi_j - \phi^\star_j\|^2_n \ge \frac18 \left(\frac{r_n}{4} - C_2 \lambda_{n, 1}\right)^2 \ge \frac{r_n^2}{2^7} \,.
\end{equation}
This yields contradiction as soon as $\tilde C^2 \ge 2^7 C$. Now we rigorize this intuition. \textcolor{black}{Define the following events: 
\begin{enumerate}
    \item $\Omega_{n, 1} = \{\sum_{j \in S_1} \|\hat \phi_j - \phi_j^\star\|_n \le C(s_1\lambda_{n, 1} + s_2\lambda_{n, 2})\}$. 
    \item $\Omega_{n, 2} = \left\{\left|\|\phi^\star_j\|_n - \|\phi_j^\star\|_2\right| \le \frac{\|\phi_j^\star\|_2^2}{2} \ \ \forall \ \ j \in S_1\right\}$.
     \item $\Omega_{n, 3} = \left\{\left|\|\phi^\star_{jk}\|_n - \|\phi_{jk}^\star\|_2\right| \le \frac{\|\phi_{jk}^\star\|_2^2}{2} \ \ \forall \ \ (j < k) \in S_2\right\}$.
\end{enumerate}
Using the proof of Theorem \ref{thm:statistical_error_lowdim} $\bbP(\Omega_{n, 1}) \to 1$ as $n \to \infty$. For $\Omega_{n, 2}$ note that: 
\begin{align*}
    \sum_{j \in S_1} \bbP\left(\left|\|\phi^\star_j\|^2_n - \|\phi_j^\star\|^2_2\right| \ge \frac{\|\phi_j^\star\|_2^2}{2}\right) & \le \sum_{j \in S_1}\frac{4\var\left(\left(\phi_j^\star(X)\right)^2\right)}{n\|\phi_j^\star\|_2^4} \\
    &\le \sum_{j \in S_1} \frac{4\bbE\left(\left(\phi_j^\star(X)\right)^4\right)}{n\|\phi_j^\star\|_2^4} \\
    & \le \sum_{j \in S_1} \frac{4B^2\bbE\left(\left(\phi_j^\star(X)\right)^2\right)}{n\|\phi_j^\star\|_2^4} \hspace{0.3in} [\text{As } \|\phi_j^\star\|_\infty \le B]\\
    & = \frac{4 s_1 B^2}{n\|\phi_j^\star\|_2^2} \le \frac{4 s_1 B^2}{nr_n^2} \,.
\end{align*}
Now by our choice of $r_n$, it is immediate that $\frac{nr_n^2}{s_1} \gtrsim \lambda_{n, 1}^2 \gtrsim \log n \to \infty$, using \eqref{eq:define_lambda}. This implies $\bbP(\Omega_{n, 2}^c) \to 0$. Similar calculation yields that $\bbP(\Omega_{n, 3}^c) \to 0$. 
Now, on the event $\Omega_{n, 1} \cap \Omega_{n, 2} \cap \Omega_{n, 3}$, we have from \eqref{eq:fixed_dim_prop} that $S_i \subseteq \hat S_i$ for $i = 1, 2$. 
}
\\\\
For the rest of the calculation, define the event $\mathcal{T}=\{S_i \subseteq \hat S_i, i=1,2\}$. 
We have shown $\mathbb{P}(\mathcal{T})=1-o(1)$. For the rest of the analysis, we assume $\mathcal{T}$ happens. Define $|\hat S_1 \cap S_1^c| = \gamma_1$ and $|\hat S_2 \cap S_2^c| = \gamma_2$. Therefore $|\hat S_1| = s_1 + \gamma_1$ and $|\hat S_2| = s_2 + \gamma_2$. Note that the values of $\gamma_1, \gamma_2$ satisfies: 
\begin{align*}
    2C_{10}\left(\gamma_1 \lambda_{n, 1}^2 + \gamma_2 \lambda_{n, 2}^2\right) & \le \frac{\lambda_{n, 1}}{2} \sum_{j \in S_1^c}\|\hat \phi_j \|_{n} + \frac{\lambda_{n, 2}}{2}\sum_{(j < k) \in S_2^c}\|\hat \phi_{jk}\|_{n} \\
    & \le 2C_{10}\left(s_1\lambda_{n, 1}^2 + s_2 \lambda_{n, 2}^2\right) 
\end{align*}
i.e. 
\begin{equation}
\label{eq:gamma_s_rel}
\gamma_1 \lambda_{n, 1}^2 + \gamma_2 \lambda_{n, 2}^2 \le s_1\lambda_{n, 1}^2 + s_2 \lambda_{n, 2}^2 \,.
\end{equation}
According to Algorithm \ref{algo:random_design_algo}, we use the second half of the data to estimate the mean function by restricting ourselves only on $\hat S_1$ and $\hat S_2$ and minimizing (unpenalized) squared error loss, i.e. our final estimate is: 
$$
\hat f^{\text{final}}(x) = \sum_{j \in \hat S_1} \hat \phi^{\fn}_j(x_j) + \sum_{(j < k) \in \hat S_2} \hat \phi_{jk}^{\fn}(x_j, x_k),
$$
where the component functions are estimated as: 
$$
\{\hat \phi_j^{\fn}\}, \{ \hat \phi_{jk}^{\fn}\} = \argmin_{\phi_j, \phi_{jk}} \frac1n \sum_i \left(Y_i - \sum_{j \in \hat S_1} \phi_j(X_{ij}) - \sum_{(j < k) \in \hat S_2} \phi_{jk}(X_{ij}, X_{ik})\right)^2 \,.
$$
Rest of the proof is similar to that of Theorem \ref{thm:statistical_error_lowdim}. By \cite[Theorem 4.8]{fan2022factor}, we know that if $f$ is an $d$-variate $\beta$ smooth function, then, there exists neural network $g$ with width $N$ and constant depth (where the constant depends on $d, \beta$) such that $$
\|f - g\|_\infty \le c N^{-\frac{2\beta}{d}} \,.
$$
for some constant $c=c(\beta, d)>0$. Here we fit univariate and bivariate functions seprately. Suppose we fit the univariate components using neural networks of width $N_1$ (and constant depth $c_1$) and bivariate components using neural networks of width $N_2$ (and constant depth $c_2$). Then to fit the additive functions, we need $|\hat S_1|$ many such univariate components and $|\hat S_2|$ many bivariate components. Therefore total number of weights (no. of active parameters) required to fit such a mean function is $W = C_1(|\hat S_1|N_1^2 + |\hat S_2|N_2^2)$. From \cite{bartlett2019nearly}, that VC-dim of such neural networks is $V_n \le C_2 W\log{W}$ (as the depth is $O(1)$). The bias-variance decomposition yields: 
$$
\|\hat f^{\fn} - f_0\|_{L_2(P_X)}^2 \le 2\left( \underbrace{\|\phi^\star - f_0\|_{L_2(P_X)}^2}_{\text{bias}} + \underbrace{\|\hat f^{\fn} - \phi^\star\|_{L_2(P_X)}^2}_{\text{variance}} \right)
$$
where $\phi^\star$ is the best approximator of $f_0$ among the class of neural network over which we are optimizing (i.e. sum of $|\hat S_1|$ many univariate networks with width $N_1$ and sum of $|\hat S_2|$ many bivariate components of width $N_2$). As $S_i \subseteq \hat S_i$ for $i \in \{1, 2\}$, we know from the proof of Theorem \ref{thm:fixed_error_sparse} (see equation \eqref{eq:define_nn_approx_linf}) that
$$
\|\phi^\star - f_0\|_{L_2(P_X)}^2 \le C_1\left(s_1 N_1^{-4\beta_1} + s_2 N_2^{-2\beta_2}\right) \,.
$$
Now to bound the variance term we use VC dimension techniques similar to the proof of Theorem \ref{thm:statistical_error_lowdim}. For any choice of component function, we can treat the overall sum $\sum_{j \in \hat S_1} \phi_j + \sum_{(j < k) \in S_2} \phi_{jk}$ as a function from the VC class with $VC$-dim $\le C_2 W\log{W}$. Therefore, from the proof of Theorem \ref{thm:statistical_error_lowdim}, we have:
\begin{align*}
\|\hat f^{\fn} - f_0\|_{L_2(P_X)}^2 & = O_p\left(\frac{V_n \log{n}}{n} + s_1 N_1^{-4\beta_1} + s_2 N_2^{-2\beta_2}\right) \\
& = O_p\left(\frac{W \log{W}\log{n}}{n} + s_1 N_1^{-4\beta_1} + s_2 N_2^{-2\beta_2}\right) \\
& = O_p\left(\frac{(s_1 + \gamma_1)N_1^2 + (s_2 + \gamma_2) N_2^2\log{\left((s_1 + \gamma_1)N_1^2 + (s_2 + \gamma_2) N_2^2\right)}\log{n}}{n} \right. \\
& \hspace{15em} \left. + s_1 N_1^{-4\beta_1} + s_2 N_2^{-2\beta_2}\right) \\
& = O_p\left((s_1 + \gamma_1)\lambda_{n, 1}^2 +  (s_2 + \gamma_2)\lambda_{n, 2}^2\right) = O_p\left(s_1\lambda_{n, 1}^2 +  s_2\lambda_{n, 2}^2\right) \,,
\end{align*}
where the second last line follows from the definition of $\lambda_{n, 1}$ and $\lambda_{n, 2}$ (see \eqref{eq:define_lambda}) and the last line follows from \eqref{eq:gamma_s_rel}. This completes the proof. 

\subsection{Proof of Theorem \ref{thm:minimax_high_dim}}
\label{sec:proof_minimax_high_dim}
Here, we extend our proof of Theorem \ref{thm:minimax_lower_bound_growing_d} to the sparse interaction model. The main technical change here is to incorporate sparsity. For that, we need slightly different definitions for $m_1, m_2$ from \eqref{eq:define_mi}. 
We define $m_1$ to be the solution of: 
\begin{align}
\label{eq:m1_def_high}
m_1^{-2\beta_1} & = c_1\left(n^{-\frac{2\beta_1}{2\beta_1 + 1}} \vee 8\frac{\log{\left(\frac{2d}{s_1} - 2\right)}}{n} \right) = \frac{c_1}{n}\left(n^{\frac{1}{2\beta_1 + 1}} \vee 8\log{\left(\frac{2d}{s_1} - 2\right)} \right)
\end{align}
and $m_2$ to be the solution of the equation: 
\begin{align}
    \label{eq:m2_def_high}
    m_2^{-2\beta_2} & = c_2\left(n^{-\frac{2\beta_2}{2\beta_2 + 2}} \vee 8\frac{\log{\left(\frac{d(d-1)}{s_2} - 2\right)}}{n}  \right) = \frac{c_2}{n}\left(n^{-\frac{2\beta_2}{2\beta_2 + 2}} \vee 8\log{\left(\frac{d(d-1)}{s_2} - 2\right)}  \right) 
\end{align}
If $m_1, m_2$ are not integers, we will take the nearest integer. As this does not affect us asymptotically, we henceforth assume them to be exact solution for the simplicity of proof. The constant $c_1, c_2 >0$ will be chosen at the end of the proof. Similar to \eqref{eq:define_mi}, define $h_i:=m_i^{-1}$, $i \in \{1, 2\}$. Consider the set $\Omega_i = \{0, 1\}^{m_i}$ for $i \in \{1, 2\}$. By Proposition \ref{prop:VG}, we can find $S_i \subset \Omega_i$ such that $|S_i| \ge 2^{m_i/8}$ and for any $\bomega, \bomega' \in S_i$, we have $\rho(\bomega, \bomega') \ge m_i/8$, where $\rho$ is the Hamming distance. As before define $\Omega = \Omega_1 \times \Omega_2$ and $S = S_1 \times S_2$. For any $\bomega \in \Omega$, we write it as $\bomega = (\bomega_1, \bomega_2)$ where $\bomega_i \in \Omega_i$. For one-dimensional components, given any $\bomega_1 \in \Omega_1$, we define a function $f^{(1)}_{\omega_1}$ as: 
$$
f^{(1)}_{\bomega_1}(x) =  \sum_{k=1}^m \bomega_{1, k}\phi_k(x) \,,
$$
where $\phi_k$ defined in \eqref{eq:define_lower_phi_lin}. Define $\cF_1 := \{f^{(1)}_{\bomega_1}: \bomega_1 \in S_1\}$ and $\Gamma_1 := |\cF_1| \ge 2^{m_1/8}$. Enumerate the functions in $\cF_1$ as $\cF_1 = \{f^{(1)}_1, \dots, f^{(1)}_{\Gamma_1}\}$. Note that for any $i \neq j$: 
\begin{align}
\label{eq:one_comp_highdim_bound}
    & \int_0^1 \left(f^{(1)}_i(x) - f^{(1)}_j(x)\right)^2 \ dx \notag \\
    &= \sum_{k=1}^m \left(\bomega_{i, k} - \bomega_{j, k}\right)^2 \int_0^1 \phi_k^2(x) \ dx \hspace{0.1in} [{\rm Since } \left\langle \phi_j , \phi_k \right\rangle_{L_2} = 0 \ \forall \ j \neq k ] \notag \\
    & = L_1^2 h^{2\beta_1 + 1} \|K\|_2^2 \rho(\bomega_i, \bomega_j)
\end{align}
For the two-dimensional components, given any $\bomega_2 \in \Omega_2$, we define $f^{(2)}_{\bomega_2}$ as: 
$$
f^{(2)}_{\bomega_2}(x, y) = \sum_{k, l=1}^m \bomega_{2, k, l}\phi_{k,l}(x, y) \,,
$$
where $\phi_{k, l}$ is defined by \eqref{eq:define_lower_phi_quad}. Define $\cF_2 := \{f^{(2)}_{\bomega_2}: \bomega_2 \in S_2\}$ and $\Gamma_2 := |\cF_2| \ge 2^{m^2_2/8}$. Enumerate the functions in $\cF_2$ as $\cF_2 = \{f^{(2)}_1, \dots, f^{(2)}_{\Gamma_2}\}$. For any $i \neq j$: 
\begin{align}
\label{eq:two_comp_highdim_bound}
    & \int_{[0, 1]^2} \left(f^{(2)}_{i}(x,y) - f^{(2)}_j(x, y)\right)^2 \ dx \ dy \notag \\
    & = \sum_{1 \le k, l \le m_2} \left(\bomega_{i, k, l} - \bomega_{j, k, l}\right)^2 \iint_0^1 \phi_{k, l}^2(x, y) \ dx \hspace{0.1in} [\text{Since }\left\langle \phi_{j, k} , \phi_{l, m} \right\rangle_{L_2} = 0 \ \forall \ (j, k) \neq (l, m)] \notag \\
    & = L_2^2 h^{2\beta_2 + 2} \|K\|_4^2 \rho(\bomega_i, \bomega_j)
\end{align}
To construct our set of sparse alternatives, define sets $U^*_1$ and $U^*_2$ as: 
\begin{align*}
    U^*_i & = \left\{u \in \{0, 1, \dots, \Gamma_1\}^{\binom{d}{2}}: \|u\|_0 = s_i\right\}, \qquad i=1,2.
\end{align*}
It is immediate that $|U^*_i| = \binom{\binom{d}{i}}{s_i}\Gamma_i^{s_i}$, $i=1,2$. Choose $U_i \subset U_i^*$ such that for any $u, v \in U_i$, we have $\rho(u, v) \ge s_i/2$. By the proof of \cite[Lemma 4]{raskutti2012minimax}, we have: 
\begin{align*}
    |U_i| & \ge \frac12 \left(\frac{\binom{d}{i}-s_i}{s_i/2}\right)^{s_i/2}\Gamma_i^{s_i/2} \,, \qquad i =1,2.
\end{align*}
For $i \le 2$, $u_i \in U_i$, define functions $f^{(i)}_{u_i}$ as: 
$$
f^{(1)}_{u_1}(X) = \sum_{j=1}^d f^{(1)}_{u_{1, j}}(X_j), \qquad f^{(2)}_{u_2}(X) = \sum_{1 \le i < j \le d} f^{(2)}_{u_{2, i, j}}(X_i, X_j)
$$
where $f^{(i)}_0 = 0$, else it is chosen from $\cF_i$. 
Now, our set of alternatives are: 
$$
f(X) = f^{(1)}_{u_1}(X) + f^{(2)}_{u_2}(X) \ \ u_1 \in U_1, u_2 \in U_2 \,.
$$
Define $\cF^{\text{sparse}}$ to be collection of all such functions $f$. Let $M = |\cF^{\text{sparse}}|$. Pick any two $f, f' \in \cF^{\text{sparse}}$. We have: 
\begin{align*}
    d^2(f, f') & = \int_{[0, 1]^d} \left(\sum_{j=1}^d (f^{(1)}_{u_{1, j}} - f^{(1)}_{u'_{1, j}})(X_j) + \sum_{1 \le i < j \le d} (f^{(2)}_{u_{2, i, j}} - f^{(2)}_{u'_{2, i, j}})(X_i, X_j)\right)^2 \ dX \\
    & = \int_{[0, 1]^d} \left(\sum_{j=1}^d (f^{(1)}_{u_{1, j}} - f^{(1)}_{u'_{1, j}})(X_j)\right)^2 \ dX + \int_{[0, 1]^d} \left(\sum_{1 \le i < j \le d} (f^{(2)}_{u_{2, i, j}} - f^{(2)}_{u'_{2, i, j}})(X_i, X_j)\right)^2 \ dX \\
    & \qquad \qquad + 2 \int_{[0, 1]^d}\left(\sum_{j=1}^d (f^{(1)}_{u_{1, j}} - f^{(1)}_{u'_{1, j}})(X_j)\right)\left(\sum_{1 \le i < j \le d} (f^{(2)}_{u_{2, i, j}} - f^{(2)}_{u'_{2, i, j}})(X_i, X_j)\right) \ dX \\
    & \triangleq T_1 + T_2 + 2T_3
\end{align*}
We now analyze each $T_i$ separately. 
\begin{align*}
    T_1 & = \int_{[0, 1]^d} \left(\sum_{j=1}^d (f^{(1)}_{u_{1, j}} - f^{(1)}_{u'_{1, j}})(X_j)\right)^2 \ dX \\
    & = \sum_{j=1}^d  \int_0^1 \left((f^{(1)}_{u_{1, j}} - f^{(1)}_{u'_{1, j}})(X_j)\right)^2 \ dX_j \hspace{0.1in} \left[{\rm Since} \int f^{(1)}_{u_1, j}(X_j) \ dX_j = 0\right] \\
    & = L_1^2 h_1^{2\beta_1 + 1} \|K\|_2^2 \sum_{j=1}^d \rho(\bomega_{u_{1, j}}, \bomega_{u'_{1, j}}) \ \mathds{1}_{u_{1, j} \neq u'_{1, j}} \hspace{0.1in} [\text{From }\eqref{eq:one_comp_highdim_bound}]\\
    T_2 & = \int_{[0, 1]^d} \left(\sum_{1 \le i < j \le d} (f^{(2)}_{u_{2, i, j}} - f^{(2)}_{u'_{2, i, j}})(X_i, X_j)\right)^2 \ dX \\
    & = \sum_{1 \le i < j \le d} \int_{[0, 1]^2} \left( (f^{(2)}_{u_{2, i, j}} - f^{(2)}_{u'_{2, i, j}})(X_i, X_j)\right)^2 \ dX_i \ dX_j \hspace{0.1in} [\text{ As marginals of }f^{(2)} \text{ are }0 ] \\
    & = L_2^2 h_2^{2\beta_2 + 2} \sum_{1 \le i < j \le d} \rho(\bomega_{u_{2, i, j}}, \bomega_{u'_{2, i, j}}) \ \mathds{1}_{u_{2, i, j} \neq u'_{2, i, j}} \hspace{0.1in} [\text{From }\eqref{eq:two_comp_highdim_bound}]
\end{align*}
Furthermore, $T_3 = 0$ as the marginals of $f^{(2)}$ are $0$. Combining the bounds of $T_i$'s, we obtain: 
\begin{align}
\label{eq:d_val_high_dim_lower}
d^2(f, f') & = L_1^2 h_1^{2\beta_1 + 1} \|K\|_2^2 \sum_{j=1}^d \rho(\bomega_{u_{1, j}}, \bomega_{u'_{1, j}}) \ \mathds{1}_{u_{1, j} \neq u'_{1, j}} \notag \\
& \qquad \qquad + L_2^2 h_2^{2\beta_2 + 2}\|K\|_2^4 \sum_{1 \le i < j \le d} \rho(\bomega_{u_{2, i, j}}, \bomega_{u'_{2, i, j}}) \ \mathds{1}_{u_{2, i, j} \neq u'_{2, i, j}} \notag \\
& \ge \frac{L_1^2 \|K\|_2^2}{8} h_1^{2\beta_1 + 1}m_1 \sum_{j=1}^d \mathds{1}_{u_{1, j} \neq u'_{1, j}} + \frac{L_2^2 \|K\|_2^4}{8}h_2^{2\beta_2 + 2}m^2_2\sum_{1 \le i < j \le d} \mathds{1}_{u_{2, i, j} \neq u'_{2, i, j}} \notag \\
& \ge \frac{L_1^2 \|K\|_2^2}{16} h_1^{2\beta_1 + 1}m_1 s_1 + \frac{L_2^2 \|K\|_2^4}{16}h_2^{2\beta_2 + 2}m^2_2s_2 \notag \\
& \ge \frac{L_1^2 \|K\|_2^2}{16} m_1^{-2\beta_1} s_1 + \frac{L_2^2 \|K\|_2^4}{16}m_2^{-2\beta_2}s_2 := 4\delta^2 \,. 
\end{align}
Similarly we can obtain an upper bound in terms of $\delta$: 
\begin{align}
    \label{eq:d_val_high_dim_upper}
    d^2(f, f') & = L_1^2 h_1^{2\beta_1 + 1} \|K\|_2^2 \sum_{j=1}^d \rho(\bomega_{u_{1, j}}, \bomega_{u'_{1, j}}) \ \mathds{1}_{u_{1, j} \neq u'_{1, j}} \notag \\
& \qquad \qquad + L_2^2 h_2^{2\beta_2 + 2}\|K\|_2^4 \sum_{1 \le i < j \le d} \rho(\bomega_{u_{2, i, j}}, \bomega_{u'_{2, i, j}}) \ \mathds{1}_{u_{2, i, j} \neq u'_{2, i, j}} \notag \\
& \le 2L_1^2\|K\|_2^2 s_1m_1h_1^{2\beta_1 + 1} + 2 L_2^2\|K\|_2^4  s_2m^2_2 h_2^{2\beta_2 + 2} \\
& = 2L_1^2\|K\|_2^2 s_1m_1^{-2\beta_1} + 2 L_2^2\|K\|_2^4  s_2 m_2^{-2\beta_2 
} = 128 \delta^2 \,.
\end{align}
Combining the bounds obtained in equation \eqref{eq:d_val_high_dim_lower} and equation \eqref{eq:d_val_high_dim_upper} we have that for any two $f, f' \in \cF^{\text{sparse}}$: 
\begin{equation}
\label{eq:d_val_high_bounds}
    4\delta^2 \le d^2(f, f') \le 128 \delta^2 \,.
\end{equation}
which further implies that $\cF^{\text{sparse}}$ is a $2\delta$ packing set of the set of all feasible functions. An application of Fano's inequality~\cite[Proof of Theorem 2.18]{mukherjee2021optimal} yields: 
\begin{equation}
    \label{eq:Fano_high}
    \inf_{\hat f} \sup_{f,P_{X}} \bbE_{P_{X}}[(\hat f(X) - f(X))^2] \ge \delta^2 \left(1 -\frac{\frac{n}{M^2} \sum_{f, f'  \in \cF^{\text{sparse}}}KL\left(\bbP_f | \bbP_{f'}\right) + \log{2}}{\log{(M - 1)}}\right)
\end{equation}
Now as the errors are normally distributed: 
$$
KL\left(\bbP_f | \bbP_{f'}\right) = \bbE\left[\left(f(X) - f'(X)\right)^2\right] = d^2(f, f')
$$
Using the bounds of equation \eqref{eq:d_val_high_bounds}, we have $KL\left(\bbP_f | \bbP_{f'}\right) \le 128\delta^2$ for all $f, f'$. Therefore, we obtain from equation \eqref{eq:Fano_high}: 
\begin{align}
    \label{eq:Fano_high_2}
    \inf_{\hat f} \sup_{f,P_{X}} \bbE_{P_{X}}[(\hat f(X) - f(X))^2] & \ge \delta^2 \left(1 -\frac{128 n\delta^2 + \log{2}}{\log{(M - 1)}}\right) \notag \\
    & \ge \delta^2 \left(1 -\frac{128 n\delta^2 + \log{2}}{2\log{M}}\right),
\end{align}
where $M=|U_1| |U_2|$. The rest of the proof is to balance $n\delta^2$ and $\log{M}$ to say we can say $\frac{128 n\delta^2 + \log{2}}{2\log{M}} \le c^*$ for some $0 < c^* < 1$. For that it is enough to show that $n\delta^2$ and $\log{M}$ has same order as we can then chosen the constants carefully so that the ratio is within $(0, 1)$. From the definition of $\delta$, we have: 
\begin{align}
\label{eq:fano_num_order}
    n\delta^2 & = \frac{L_1^2 \|K\|_2^2}{64} nm_1^{-2\beta_1} s_1 + \frac{L_2^2 \|K\|_2^4}{64}nm_2^{-2\beta_2}s_2  \notag\\
    & = \frac{L_1^2 \|K\|_2^2}{4} nm_1^{-2\beta_1} \frac{s_1}{16} + \frac{L_2^2 \|K\|_2^4}{4}nm_2^{-2\beta_2}\frac{s_2}{16} \notag \\
    & = \frac{L_1^2 c_1\|K\|_2^2}{4}\frac{s_1}{16}\left[n^{\frac{1}{2\beta_1 + 1}} \vee 8\log{\left(\frac{2d}{s_1} - 2\right)} \right] \notag\\
    & \qquad \qquad \qquad + \frac{L_2^2c_2 \|K\|_2^4}{4}\frac{s_2}{16} \left[(n^{\frac{1}{2\beta_2 + 2}} \vee 8\log{\left(\frac{d(d-1)}{s_2} - 2\right)}  \right]
\end{align}
On the other hand,by our construction, we have: 
\begin{align}
\label{eq:fano_denom}
\log{M} & = \log{|U_1|} + \log{|U_2|} \notag \\
& \ge 2\log{\frac12} + \frac{s_1}{2} \log{\left(\frac{d-s_1}{s_1/2}\right)} + \frac{s_1}{2} \log{\Gamma_1} + \frac{s_2}{2} \log{\left(\frac{d(d-1)/2 - s_2}{s_2/2}\right)} + \frac{s_2}{2} \log{\Gamma_2} \notag \\
& = 2\log{\frac12} + \frac{s_1}{2} \log{\left(\frac{2d}{s_1} - 2\right)} + \frac{s_1}{2} \log{\Gamma_1} + \frac{s_2}{2} \log{\left(\frac{d(d-1)}{s_2} - 2\right)} + \frac{s_2}{2} \log{\Gamma_2} \notag \\
& = 2\log{\frac12} + \frac{s_1}{2} \log{\left(\frac{2d}{s_1} - 2\right)} + \frac{s_1 m_1}{16}  + \frac{s_2}{2} \log{\left(\frac{d(d-1)}{s_2} - 2\right)} + \frac{s_2 m_2^2}{16}  \notag \\
& = 2\log{\frac12} + \frac{s_1}{16}\left[m_1 + 8\log{\left(\frac{2d}{s_1} - 2\right)}\right] + 
\frac{s_2}{16}\left[m_2^2 + 8\log{\left(\frac{d(d-1)}{s_2} - 2\right)}\right] \notag \\
& \triangleq 2\log{\frac12} + a_n + b_n \,.
\end{align}
We now consider four cases: 
\begin{framed}
{\bf Case 1: $n^{1/(2\beta_1 + 1)} \ge 8\log{\left(\frac{2d}{s_1} - 2\right)}$ and $n^{\frac{1}{2\beta_2 + 2}} \ge 8\log{\left(\frac{d(d-1)}{s_2} - 2\right)}$}
\end{framed}
Using \eqref{eq:m1_def_high}, \eqref{eq:m2_def_high}, we have $m_1 = n^{1/(2\beta_1 + 1)}$ and $m_2 = n^{1/(2\beta_2 + 2)}$. Therefore, equations \eqref{eq:fano_num_order} and \eqref{eq:fano_denom} are simplified to:
\begin{align*}
    n\delta^2 & = \frac{L_1^2 c_1\|K\|_2^2}{4}\frac{s_1}{16}n^{\frac{1}{2\beta_1 + 1}} +\frac{L_2^2c_2 \|K\|_2^4}{4}\frac{s_2}{16} n^{-\frac{2\beta_2}{2\beta_2 + 2}} \\
    & \le \frac{L_1^2 c_1\|K\|_2^2}{4}\frac{s_1}{16}\left[m_1 + 8\log{\left(\frac{2d}{s_1} - 2\right)}\right]  +\frac{L_2^2c_2 \|K\|_2^4}{4}\frac{s_2}{16}\left[m_2^2 + 8\log{\left(\frac{d(d-1)}{s_2} - 2\right)}\right] \\
    & \le \frac{L_1^2 c_1\|K\|_2^2}{4}a_n +\frac{L_2^2c_2 \|K\|_2^4}{4}b_n
\end{align*}

\begin{framed}
{\bf Case 2: $n^{1/(2\beta_1 + 1)} \ge 8\log{\left(\frac{2d}{s_1} - 2\right)}$ and $n^{\frac{1}{2\beta_2 + 2}} < 8\log{\left(\frac{d(d-1)}{s_2} - 2\right)}$}
\end{framed}
Here, $m_1 = n^{1/(2\beta_1 + 1)}$ and $m_2 =\left(8\log{\left(\frac{d(d-1)}{s_2} - 2\right)}/n\right)^{-1/2\beta_2}$ yielding
\begin{align*}
    n\delta^2 & = \frac{L_1^2 c_1\|K\|_2^2}{4}\frac{s_1}{16}n^{\frac{1}{2\beta_1 + 1}} +\frac{L_2^2c_2 \|K\|_2^4}{4}\frac{s_2}{16} 8\log{\left(\frac{d(d-1)}{s_2} - 2\right)} \\
    & \le \frac{L_1^2 c_1\|K\|_2^2}{4}\frac{s_1}{16}\left[m_1 + 8\log{\left(\frac{2d}{s_1} - 2\right)}\right]  +\frac{L_2^2c_2 \|K\|_2^4}{4}\frac{s_2}{16}\left[m_2^2 + 8\log{\left(\frac{d(d-1)}{s_2} - 2\right)}\right] \\
    & \le \frac{L_1^2 c_1\|K\|_2^2}{4}a_n +\frac{L_2^2c_2 \|K\|_2^4}{4}b_n
\end{align*}

\begin{framed}
{\bf Case 3: $n^{1/(2\beta_1 + 1)} < 8\log{\left(\frac{2d}{s_1} - 2\right)}$ and $n^{\frac{1}{2\beta_2 + 2}} \ge 8\log{\left(\frac{d(d-1)}{s_2} - 2\right)}$}
\end{framed}
Here $m_1 = \left(8\log{\left(\frac{2d}{s_1} - 2\right)}/n\right)^{-1/2\beta_1}$ and $m_2 =n^{\frac{1}{2\beta_2 + 2}}$ yielding
\begin{align*}
    n\delta^2 & = \frac{L_1^2 c_1\|K\|_2^2}{4}\frac{s_1}{16}8\log{\left(\frac{2d}{s_1} - 2\right)} +\frac{L_2^2c_2 \|K\|_2^4}{4}\frac{s_2}{16} n^{\frac{2}{2\beta_2 + 2}}\\
    & \le \frac{L_1^2 c_1\|K\|_2^2}{4}\frac{s_1}{16}\left[m_1 + 8\log{\left(\frac{2d}{s_1} - 2\right)}\right]  +\frac{L_2^2c_2 \|K\|_2^4}{4}\frac{s_2}{16}\left[m_2^2 + 8\log{\left(\frac{d(d-1)}{s_2} - 2\right)}\right] \\
    & \le \frac{L_1^2 c_1\|K\|_2^2}{4}a_n +\frac{L_2^2c_2 \|K\|_2^4}{4}b_n.
\end{align*}

\begin{framed}
{\bf Case 4: $n^{1/(2\beta_1 + 1)} < 8\log{\left(\frac{2d}{s_1} - 2\right)}$ and $n^{\frac{1}{2\beta_2 + 2}} < 8\log{\left(\frac{d(d-1)}{s_2} - 2\right)}$}
\end{framed}
Here $m_1 = \left(8\log{\left(\frac{2d}{s_1} - 2\right)}/n\right)^{-1/2\beta_1}$ and $m_2 =\left(8\log{\left(\frac{d(d-1)}{s_2} - 2\right)}/n\right)^{-1/2\beta_2}$ yielding
\begin{align*}
    n\delta^2 & = \frac{L_1^2 c_1\|K\|_2^2}{4}\frac{s_1}{16}8\log{\left(\frac{2d}{s_1} - 2\right)} +\frac{L_2^2c_2 \|K\|_2^4}{4}\frac{s_2}{16} 8\log{\left(\frac{d(d-1)}{s_2} - 2\right)}\\
    & \le \frac{L_1^2 c_1\|K\|_2^2}{4}\frac{s_1}{16}\left[m_1 + 8\log{\left(\frac{2d}{s_1} - 2\right)}\right]  +\frac{L_2^2c_2 \|K\|_2^4}{4}\frac{s_2}{16}\left[m_2^2 + 8\log{\left(\frac{d(d-1)}{s_2} - 2\right)}\right] \\
    & \le \frac{L_1^2 c_1\|K\|_2^2}{4}a_n +\frac{L_2^2c_2 \|K\|_2^4}{4}b_n.
\end{align*}
\\
Hence, in all four cases, the numerator and the denominator has same order. Hence,
$$
n\delta^2 \le \frac{L_1^2 c_1\|K\|_2^2}{4}a_n +\frac{L_2^2c_2 \|K\|_2^4}{4}b_n \,.
$$
Putting this in equation \eqref{eq:Fano_high_2} we obtain: 
\begin{align}
    \inf_{\hat f} \sup_{f,P_{X}} \bbE_{P_{X}}[(\hat f(X) - f(X))^2] & \ge \delta^2 \left(1 -\frac{32 \left(L_1^2 c_1\|K\|_2^2 \ a_n + L_2^2c_2 \|K\|_2^4 \ b_n\right) + \log{2}}{2(2\log{\frac12} + a_n + b_n)}\right) \notag \\
    & \ge \delta^2 \left(1 -\frac{33 \left(L_1^2 c_1\|K\|_2^2 \ a_n + L_2^2c_2 \|K\|_2^4 \ b_n\right)}{3(a_n + b_n)}\right) \notag \\
    & \ge \delta^2 \left(1 -\frac{11 \left(L_1^2 c_1\|K\|_2^2 \ a_n + L_2^2c_2 \|K\|_2^4 \ b_n\right)}{(a_n + b_n)}\right). \notag \\
\end{align}
Now choose $c_1, c_2$ such that $(L_1^2 c_1\|K\|_2^2) \vee (L_2^2c_2 \|K\|_2^4) \le 1/22$. Then we have: 
\begin{align*}
     & \inf_{\hat f} \sup_{f,P_{X}} \bbE_{P_{X}}[(\hat f(X) - f(X))^2] \\
     & \ge \frac12 \delta^2 \ge \frac{L_1^2 \|K\|_2^2}{128} m_1^{-2\beta_1} s_1 + \frac{L_2^2 \|K\|_2^4}{128}m_2^{-2\beta_2}s_2 \\
     & = \frac{L_1^2 \|K\|_2^2c_1}{128} s_1\left(n^{-\frac{2\beta_1}{2\beta_1 + 1}} \vee 8\frac{\log{\left(\frac{2d}{s_1} - 2\right)}}{n} \right) + \frac{L_2^2 \|K\|_2^4 c_2}{128}s_2\left(n^{-\frac{2\beta_2}{2\beta_2 + 2}} \vee 8\frac{\log{\left(\frac{d(d-1)}{s_2} - 2\right)}}{n}  \right) .
\end{align*}
This completes the proof. 

\end{document}